\documentclass{article}
\usepackage{verbatim,amsmath,amscd,amssymb,amsthm}
\usepackage[all]{xy}
\begin{document}
\font\germ=eufm10
\def\ssl{\hbox{\germ sl}}
\def\slh{\widehat{\ssl_2}}
\def\ge{\hbox{\germ g}}
\def\aaa{@}
\title{\large\bf GEOMETRIC CRYSTALS ON UNIPOTENT GROUPS
 AND \\ GENERALIZED YOUNG TABLEAUX }

\vskip7pt
\author{\large N\textsc{akashima} Toshiki
\thanks
{supported in part by JSPS Grants in Aid for 
Scientific Research}
\\
Department of Mathematics, 
\\
Sophia University, Tokyo 102-8554, Japan\\
e-mail: toshiki@mm.sophia.ac.jp
}
\date{}
\maketitle
\begin{abstract}
We define geometric/unipotent crystal structure on 
unipotent subgroups of semi-simple algebraic groups.
We shall show that in $A_n$-case, their ultra-discretizations 
coincide with crystals obtained by generalizing 
Young tableaux.
\end{abstract}


\renewcommand{\labelenumi}{$($\roman{enumi}$)$}
\renewcommand{\labelenumii}{$(${\rm \alph{enumii}}$)$}
\font\germ=eufm10
\def\m@th{\mathsurround=0pt}

\def\m@th{\mathsurround=0pt}

\def\fsquare(#1,#2){
\hbox{\vrule$\hskip-0.4pt\vcenter to #1{\normalbaselines\m@th
\hrule\vfil\hbox to #1{\hfill$\scriptstyle #2$\hfill}\vfil\hrule}$\hskip-0.4pt
\vrule}}

\def\addsquare(#1,#2){\hbox{$
	\dimen1=#1 \advance\dimen1 by -0.8pt
	\vcenter to #1{\hrule height0.4pt depth0.0pt%
	\hbox to #1{%
	\vbox to \dimen1{\vss%
	\hbox to \dimen1{\hss$\scriptstyle~#2~$\hss}%
	\vss}%
	\vrule width0.4pt}%
	\hrule height0.4pt depth0.0pt}$}}

\def\Fsquare(#1,#2){
\hbox{\vrule$\hskip-0.4pt\vcenter to #1{\normalbaselines\m@th
\hrule\vfil\hbox to #1{\hfill$#2$\hfill}\vfil\hrule}$\hskip-0.4pt
\vrule}}

\def\Addsquare(#1,#2){\hbox{$
	\dimen1=#1 \advance\dimen1 by -0.8pt
	\vcenter to #1{\hrule height0.4pt depth0.0pt%
	\hbox to #1{%
	\vbox to \dimen1{\vss%
	\hbox to \dimen1{\hss$~#2~$\hss}%
	\vss}%
	\vrule width0.4pt}%
	\hrule height0.4pt depth0.0pt}$}}

\def\hfourbox(#1,#2,#3,#4){%
	\fsquare(0.3cm,#1)\addsquare(0.3cm,#2)\addsquare(0.3cm,#3)\addsquare(0.3cm,#4)}

\def\Hfourbox(#1,#2,#3,#4){%
	\Fsquare(0.4cm,#1)\Addsquare(0.4cm,#2)\Addsquare(0.4cm,#3)\Addsquare(0.4cm,#4)}

\def\HHfourbox(#1,#2,#3,#4){%
	\Fsquare(0.8cm,#1)\Addsquare(0.8cm,#2)\Addsquare(0.8cm,#3)\Addsquare(0.8cm,#4)}

\def\hthreebox(#1,#2,#3){%
	\fsquare(0.3cm,#1)\addsquare(0.3cm,#2)\addsquare(0.3cm,#3)}

\def\htwobox(#1,#2){%
	\fsquare(0.3cm,#1)\addsquare(0.3cm,#2)}

\def\vfourbox(#1,#2,#3,#4){%
	\normalbaselines\m@th\offinterlineskip
	\vcenter{\hbox{\fsquare(0.3cm,#1)}
	      \vskip-0.4pt
	      \hbox{\fsquare(0.3cm,#2)}	
	      \vskip-0.4pt
	      \hbox{\fsquare(0.3cm,#3)}	
	      \vskip-0.4pt
	      \hbox{\fsquare(0.3cm,#4)}}}

\def\VVfourbox(#1,#2,#3,#4){%
	\normalbaselines\m@th\offinterlineskip
	\vcenter{\hbox{\Fsquare(0.8cm,#1)}
	      \vskip-0.4pt
	      \hbox{\Fsquare(0.8cm,#2)}	
	      \vskip-0.4pt
	      \hbox{\Fsquare(0.8cm,#3)}	
	      \vskip-0.4pt
	      \hbox{\Fsquare(0.8cm,#4)}}}

\def\Vfourbox(#1,#2,#3,#4){%
	\normalbaselines\m@th\offinterlineskip
	\vcenter{\hbox{\Fsquare(0.4cm,#1)}
	      \vskip-0.4pt
	      \hbox{\Fsquare(0.4cm,#2)}	
	      \vskip-0.4pt
	      \hbox{\Fsquare(0.4cm,#3)}	
	      \vskip-0.4pt
	      \hbox{\Fsquare(0.4cm,#4)}}}

\def\vthreebox(#1,#2,#3){%
	\normalbaselines\m@th\offinterlineskip
	\vcenter{\hbox{\fsquare(0.3cm,#1)}
	      \vskip-0.4pt
	      \hbox{\fsquare(0.3cm,#2)}	
	      \vskip-0.4pt
	      \hbox{\fsquare(0.3cm,#3)}}}

\def\vtwobox(#1,#2){%
	\normalbaselines\m@th\offinterlineskip
	\vcenter{\hbox{\fsquare(0.3cm,#1)}
	      \vskip-0.4pt
	      \hbox{\fsquare(0.3cm,#2)}}}

\def\Hthreebox(#1,#2,#3){%
	\Fsquare(0.4cm,#1)\Addsquare(0.4cm,#2)\Addsquare(0.4cm,#3)}

\def\HHthreebox(#1,#2,#3){%
	\Fsquare(0.8cm,#1)\Addsquare(0.8cm,#2)\Addsquare(0.8cm,#3)}

\def\Htwobox(#1,#2){%
	\Fsquare(0.4cm,#1)\Addsquare(0.4cm,#2)}

\def\H6twobox(#1,#2){%
	\Fsquare(0.6cm,#1)\Addsquare(0.6cm,#2)}

\def\HHtwobox(#1,#2){%
	\Fsquare(0.8cm,#1)\Addsquare(0.8cm,#2)}

\def\Vthreebox(#1,#2,#3){%
	\normalbaselines\m@th\offinterlineskip
	\vcenter{\hbox{\Fsquare(0.4cm,#1)}
	      \vskip-0.4pt
	      \hbox{\Fsquare(0.4cm,#2)}	
	      \vskip-0.4pt
	      \hbox{\Fsquare(0.4cm,#3)}}}

\def\Vtwobox(#1,#2){%
	\normalbaselines\m@th\offinterlineskip
	\vcenter{\hbox{\Fsquare(0.4cm,#1)}
	      \vskip-0.4pt
	      \hbox{\Fsquare(0.4cm,#2)}}}

\def\twoone(#1,#2,#3){%
	\normalbaselines\m@th\offinterlineskip
	\vcenter{\hbox{\htwobox({#1},{#2})}
	      \vskip-0.4pt
	      \hbox{\fsquare(0.3cm,#3)}}}

\def\Twoone(#1,#2,#3){%
	\normalbaselines\m@th\offinterlineskip
	\vcenter{\hbox{\Htwobox({#1},{#2})}
	      \vskip-0.4pt
	      \hbox{\Fsquare(0.4cm,#3)}}}

\def\TTwoone(#1,#2,#3){%
	\normalbaselines\m@th\offinterlineskip
	\vcenter{\hbox{\H6twobox({#1},{#2})}
	      \vskip-0.4pt
	      \hbox{\Fsquare(0.6cm,#3)}}}

\def\threeone(#1,#2,#3,#4){%
	\normalbaselines\m@th\offinterlineskip
	\vcenter{\hbox{\hthreebox({#1},{#2},{#3})}
	      \vskip-0.4pt
	      \hbox{\fsquare(0.3cm,#4)}}}

\def\Threeone(#1,#2,#3,#4){%
	\normalbaselines\m@th\offinterlineskip
	\vcenter{\hbox{\Hthreebox({#1},{#2},{#3})}
	      \vskip-0.4pt
	      \hbox{\Fsquare(0.4cm,#4)}}}

\def\Threetwo(#1,#2,#3,#4,#5){%
	\normalbaselines\m@th\offinterlineskip
	\vcenter{\hbox{\Hthreebox({#1},{#2},{#3})}
	      \vskip-0.4pt
	      \hbox{\Htwobox({#4},{#5})}}}

\def\threetwo(#1,#2,#3,#4,#5){%
	\normalbaselines\m@th\offinterlineskip
	\vcenter{\hbox{\hthreebox({#1},{#2},{#3})}
	      \vskip-0.4pt
	      \hbox{\htwobox({#4},{#5})}}}

\def\twotwo(#1,#2,#3,#4){%
	\normalbaselines\m@th\offinterlineskip
	\vcenter{\hbox{\htwobox({#1},{#2})}
	      \vskip-0.4pt
	      \hbox{\htwobox({#3},{#4})}}}

\def\Twotwo(#1,#2,#3,#4){%
	\normalbaselines\m@th\offinterlineskip
	\vcenter{\hbox{\Htwobox({#1},{#2})}
	      \vskip-0.4pt
	      \hbox{\Htwobox({#3},{#4})}}}

\def\TTwotwo(#1,#2,#3,#4){%
	\normalbaselines\m@th\offinterlineskip
	\vcenter{\hbox{\H6twobox({#1},{#2})}
	      \vskip-0.4pt
	      \hbox{\H6twobox({#3},{#4})}}}

\def\twooneone(#1,#2,#3,#4){%
	\normalbaselines\m@th\offinterlineskip
	\vcenter{\hbox{\htwobox({#1},{#2})}
	      \vskip-0.4pt
	      \hbox{\fsquare(0.3cm,#3)}
	      \vskip-0.4pt
	      \hbox{\fsquare(0.3cm,#4)}}}

\def\Twooneone(#1,#2,#3,#4){%
	\normalbaselines\m@th\offinterlineskip
	\vcenter{\hbox{\Htwobox({#1},{#2})}
	      \vskip-0.4pt
	      \hbox{\Fsquare(0.4cm,#3)}
	      \vskip-0.4pt
	      \hbox{\Fsquare(0.4cm,#4)}}}

\def\Twotwoone(#1,#2,#3,#4,#5){%
	\normalbaselines\m@th\offinterlineskip
	\vcenter{\hbox{\Htwobox({#1},{#2})}
	      \vskip-0.4pt
	      \hbox{\Htwobox({#3},{#4})}
              \vskip-0.4pt
	      \hbox{\Fsquare(0.4cm,#5)}}}

\def\twotwoone(#1,#2,#3,#4,#5){%
	\normalbaselines\m@th\offinterlineskip
	\vcenter{\hbox{\htwobox({#1},{#2})}
	      \vskip-0.4pt
	      \hbox{\htwobox({#3},{#4})}
              \vskip-0.4pt
	      \hbox{\fsquare(0.3cm,#5)}}}

\def\twotwotwo(#1,#2,#3,#4,#5,#6){%
	\normalbaselines\m@th\offinterlineskip
	\vcenter{\hbox{\htwobox({#1},{#2})}
	      \vskip-0.4pt
	      \hbox{\htwobox({#3},{#4})}
              \vskip-0.4pt
	      \hbox{\htwobox({#5},{#6})}}}

\def\threetwoone(#1,#2,#3,#4,#5,#6){%
	\normalbaselines\m@th\offinterlineskip
	\vcenter{\hbox{\hthreebox({#1},{#2},{#3})}
	      \vskip-0.4pt
	      \hbox{\htwobox({#4},{#5})}
              \vskip-0.4pt
	      \hbox{\fsquare(0.3cm,#6)}}}

\def\TThreetwoone(#1,#2,#3,#4,#5,#6){%
	\normalbaselines\m@th\offinterlineskip
	\vcenter{\hbox{\HHthreebox({#1},{#2},{#3})}
	      \vskip-0.4pt
	      \hbox{\HHtwobox({#4},{#5})}
              \vskip-0.4pt
	      \hbox{\fsquare(0.8cm,#6)}}}

\def\Threeoneone(#1,#2,#3,#4,#5){%
	\normalbaselines\m@th\offinterlineskip
	\vcenter{\hbox{\Hthreebox({#1},{#2},{#3})}
	      \vskip-0.4pt
	      \hbox{\Fsquare(0.4cm,#4)}
              \vskip-0.4pt
	      \hbox{\Fsquare(0.4cm,#5)}}}

\def\threeoneone(#1,#2,#3,#4,#5){%
	\normalbaselines\m@th\offinterlineskip
	\vcenter{\hbox{\hthreebox({#1},{#2},{#3})}
	      \vskip-0.4pt
	      \hbox{\fsquare(0.3cm,#4)}
              \vskip-0.4pt
	      \hbox{\fsquare(0.3cm,#5)}}}

\def\FFourtwoone(#1,#2,#3,#4,#5,#6,#7){%
	\normalbaselines\m@th\offinterlineskip
	\vcenter{\hbox{\HHfourbox({#1},{#2},{#3},{#4})}
	      \vskip-0.4pt
	      \hbox{\HHtwobox({#5},{#6})}
              \vskip-0.4pt
	      \hbox{\fsquare(0.8cm,#7)}}}

\def\a{\fsquare(0.3cm){1}\addsquare(0.3cm)(2)\addsquare(0.3cm)(3)}

\def\b{\hbox{%
	\normalbaselines\m@th\offinterlineskip
	\vcenter{\hbox{\fsquare(0.3cm){2}}\vskip-0.4pt\hbox{\fsquare(0.3cm){2}}}}}

\def\c{\hbox{\normalbaselines\m@th\offinterlineskip%
	\vcenter{\hbox{\a}\vskip-0.4pt\hbox{\b}}}}


\dimen1=0.4cm\advance\dimen1 by -0.8pt
\def\ffsquare#1{%
	\fsquare(0.4cm,\hbox{#1})}

\def\naga{%
	\hbox{$\vcenter to 0.4cm{\normalbaselines\m@th
	\hrule\vfil\hbox to 1.2cm{\hfill$\cdots$\hfill}\vfil\hrule}$}}

\def\vnaga{\normalbaselines\m@th\baselineskip0pt\offinterlineskip%
	\vrule\vbox to 1.2cm{\vskip7pt\hbox to \dimen1{$\hfil\vdots\hfil$}\vfil}\vrule}

\def\dhbox{\hbox{$\ffsquare 1 \naga \ffsquare N$}}

\def\dvbox{\hbox{\normalbaselines\m@th\baselineskip0pt\offinterlineskip\vbox{%
	  \hbox{$\ffsquare 1$}\vskip-0.4pt\hbox{$\vnaga$}\vskip-0.4pt\hbox{$\ffsquare N$}}}}

\def\sq(#1){\fsquare(0.4cm,#1)}
\def\Sq(#1){\fsquare(0.5cm,#1)}
\def\SSq(#1){\fsquare(0.9cm,#1)}
\def\Sqj{\Sq(j)}
\def\Sqjb{\Sq(\ovl j)}
\def\Sqjm{\Sq(\scriptstyle{j-1})}
\def\Sqjp{\Sq(\scriptstyle{j+1})}
\def\Sqjmb{\Sq(\scriptstyle{\ovl{j-1}})}
\def\SSqj{\SSq(j)}
\def\Sqjpb{\Sq(\scriptstyle{\ovl{j+1}})}
\def\SSqjb{\SSq(\ovl j)}
\def\SSqjm{\SSq(\scriptstyle{j-1})}
\def\SSqjp{\SSq(\scriptstyle{j+1})}
\def\SSqjmb{\SSq(\scriptstyle{\ovl{j-1}})}
\def\SSqjpb{\SSq(\scriptstyle{\ovl{j+1}})}

\def\mapright#1{\smash{\mathop{\longrightarrow}\limits^{#1}}}
\def\map#1{\smash{\mathop{\longmapsto}\limits^{#1}}}

\def\m@th{\mathsurround=0pt}

\def\black-box{
	\normalbaselines\m@th
	$\hskip-0.3pt\vrule depth 2.1pt height 7.1pt width 8pt$
	}

\def\blackbox{
	\normalbaselines\m@th
	$\hskip-0.2pt\vrule depth 2.1pt height 6.5pt width 8.6pt$
	}

\def\edge{
	\normalbaselines\m@th
	\vskip-3.4pt$\vrule depth 0.4pt height 0.1pt width 9.5pt$
	}

\def\Blackbox{
	\normalbaselines\m@th
	$\hskip-0.2pt\vrule depth 2.1pt height 8.5pt width 11.1pt$
	}

\def\Blackhalfbox{
	\normalbaselines\m@th
	$\hskip-0.2pt\vrule depth 2.1pt height 8.5pt width 5pt$
	}

\def\twoblackbox{
	\normalbaselines\m@th
	$\hskip-0.3pt\vrule depth 2.1pt height 6.5pt width 16.8pt$
	}

\def\halfoneblackbox{
	\normalbaselines\m@th
	$\hskip-0.2pt\vrule depth 2.1pt height 8.3pt width 16.3pt$
	}

\def\Fsquare(#1,#2){
\hbox{\vrule$\hskip-0.4pt\vcenter to #1{\normalbaselines\m@th
\hrule\vfil\hbox to #1{\hfill$\scriptstyle #2$\hfill}\vfil\hrule}$\hskip-0.4pt
\vrule}}

\def\halfsquare(#1,#2,#3){
\hbox{\vrule$\hskip-0.4pt\vcenter to #1{\normalbaselines\m@th
\hrule\vfil\hbox to #2{\hfill$\scriptstyle #3$\hfill}\vfil\hrule}$\hskip-1pt
\vrule}}

\def\addsquare(#1,#2){\hbox{$
	\dimen1=#1 \advance\dimen1 by -0.8pt
	\vcenter to #1{\hrule height0.4pt depth0.0pt%
	\hbox to #1{%
        \vbox to \dimen1{\vss%
	\hbox to \dimen1{\hss$\scriptstyle~#2~$\hss}%
	\vss}%
	\vrule width0.4pt}%
	\hrule height0.4pt depth0.0pt}$}}

\def\Addsquare(#1,#2){\hbox{$
	\dimen1=#1 \advance\dimen1 by -0.8pt
	\vcenter to #1{\hrule height0.4pt depth0.0pt%
	\hbox to #1{%
	\vbox to \dimen1{\vss%
	\hbox to \dimen1{\hss$\scriptstyle~#2~$\hss}%
	\vss}%
	\vrule width0.4pt}%
	\hrule height0.4pt depth0.0pt}$}}

\def\fsquare(#1,#2){
\hbox{\vrule$\hskip-0.4pt\vcenter to #1{\normalbaselines\m@th
\hrule\vfil\hbox to #1{\hfill$\scriptstyle #2$\hfill}\vfil\hrule}$\hskip-0.4pt
\vrule}}

\def\seudosquare(#1,#2,#3){
\hbox{$\hskip-0.4pt\vcenter to #1{\normalbaselines\m@th
\hrule\vfil\hbox to #2{$\hfill\scriptstyle #3\hfill$}\vfil\hrule}$\hskip-0.4pt
\vrule}}

\def\topseudosquare(#1,#2,#3){
\hbox{$\hskip-0.4pt\vtop to #1{\normalbaselines\m@th
\hrule\vfil\hbox to #2{$\hfill\scriptstyle #3\hfill$}\vfil\hrule}$\hskip-0.4pt
\vrule}}

\def\tophalfsquare(#1,#2,#3){
\hbox{\vrule$\hskip-0.4pt\vtop to #1{\normalbaselines\m@th
\hrule\vfil\hbox to #2{\hfill$\scriptstyle #3$\hfill}\vfil\hrule}$\hskip-1pt
\vrule}}

\def\Young(#1,#2,#3,#4,#5){
\vcenter{\hbox{
        $\tophalfsquare(2.0cm,0.4cm,#1)
         \topseudosquare(1.7cm,0.4cm,#2)
         \topseudosquare(1.4cm,0.4cm,#3)
         \topseudosquare(1.1cm,0.4cm,#4)
         \topseudosquare(0.8cm,0.4cm,#5)$}}}

\def\HHighest{%
	\normalbaselines\m@th\offinterlineskip
\vcenter{\hbox{\halfsquare(0.8cm,8.4cm,111\qq\cd\cd\cd\cd\qq 1111)}
	      \vskip-0.4pt
\hbox{\halfsquare(0.8cm,7.0cm,222\q\q\cd\cd\q\q 222)}
              \vskip-0.4pt
    \hbox{\halfsquare(0.8cm,5.4cm,\cd\cd\cd)}
              \vskip-0.4pt
\hbox{\halfsquare(0.8cm,3.8cm,{\syl nnn\,\cd\cd\, nnn})}
              \vskip-0.4pt
 \hbox{\halfsquare(0.8cm,3.2cm,{\syl \cd\cd})}}}

\def\topspinbox(#1,#2){%
	\normalbaselines\m@th\offinterlineskip
	\vtop{\hbox{\tophalfsquare(0.4cm,2.8mm,#1)}
	      \vskip-0.4pt
	      \hbox{\tophalfsquare(1.4cm,2.8mm,\vdots)}	
	      \vskip-0.4pt
	      \hbox{\tophalfsquare(0.4cm,2.8mm,#2)}}}

\def\spinYoung(#1,#2,#3,#4,#5){
\vcenter{\hbox{
        $\tophalfsquare(2.1cm,0.3cm,#1)
         \topseudosquare(1.7cm,0.4cm,#2)
         \topseudosquare(1.4cm,0.4cm,#3)
         \topseudosquare(1.1cm,0.4cm,#4)
         \topseudosquare(0.8cm,0.4cm,#5)$}}}

\def\bigYoung(#1,#2,#3,#4,#5){
\vcenter{\hbox{
        $\tophalfsquare(2.2cm,0.5cm,#1)
         \topseudosquare(1.9cm,0.5cm,#2)
         \topseudosquare(1.6cm,0.5cm,#3)
         \topseudosquare(1.3cm,0.5cm,#4)
         \topseudosquare(1.1cm,0.5cm,#5)$}}}

\def\YYoung(#1,#2,#3,#4,#5){
\vcenter{\hbox{
        $\tophalfsquare(3.2cm,0.7cm,#1)
         \topseudosquare(2.7cm,0.7cm,#2)
         \topseudosquare(2.2cm,0.7cm,#3)
         \topseudosquare(1.8cm,0.7cm,#4)
         \topseudosquare(1.6cm,0.7cm,#5)$}}}

\def\Longtwobox{%
	$\halfsquare(0.9cm,0.4cm,\vdots)
        \seudosquare(0.9cm,0.4cm,\vdots)$}

\def\spinlongtwobox{%
	$\halfsquare(1.3cm,0.3cm,\vdots)
        \seudosquare(1.3cm,0.4cm,\vdots)$}

\def\spintwobox(#1,#2){%
	$\halfsquare(0.4cm,0.3cm,#1)
         \seudosquare(0.4cm,0.4cm,#2)$}

\def\spinonebox(#1){%
	$\halfsquare(0.4cm,0.3cm,#1)
         \seudosquare(0.4cm,0cm,{})$}

\def\standardonebox(#1){%
	$\halfsquare(0.4cm,0.4cm,#1)
         \seudosquare(0.4cm,0cm,{})$}

\def\standardhalfbox(#1){%
	$\halfsquare(0.6cm,0.4cm,#1)
         \seudosquare(0.6cm,0cm,{})$}

\def\standardtwobox(#1,#2){%
	$\halfsquare(0.4cm,0.4cm,#1)
         \seudosquare(0.4cm,0.4cm,#2)$}

\def\Longbox(#1,#2){%
	\normalbaselines\m@th\offinterlineskip
	\vcenter{\hbox{\halfsquare(0.4cm,0.4cm,#1)}
	      \vskip-0.4pt
	      \hbox{\halfsquare(1.2cm,0.4cm,\vdots)}	
	      \vskip-0.4pt
	      \hbox{\halfsquare(0.4cm,0.4cm,#2)}}}

\def\LLongbox(#1,#2){%
	\normalbaselines\m@th\offinterlineskip
	\vcenter{\hbox{\halfsquare(0.8cm,0.8cm,#1)}
	      \vskip-0.4pt
	      \hbox{\halfsquare(1.5cm,0.8cm,\vdots)}	
	      \vskip-0.4pt
	      \hbox{\halfsquare(0.8cm,0.8cm,#2)}}}

\def\spinbox(#1,#2){%
	\normalbaselines\m@th\offinterlineskip
	\vcenter{\hbox{\halfsquare(4.2mm,2.8mm,#1)}
	      \vskip-0.4pt
	      \hbox{\halfsquare(1.4cm,2.8mm,\vdots)}	
	      \vskip-0.4pt
	      \hbox{\halfsquare(4.2mm,2.8mm,#2)}}}

\def\spinthreebox(#1,#2,#3){%
	\normalbaselines\m@th\offinterlineskip
	\vcenter{\hbox{\halfsquare(0.4cm,0.2cm,#1)}
	      \vskip-0.4pt
	      \hbox{\halfsquare(0.4cm,0.2cm,#2)}	
	      \vskip-0.4pt
	      \hbox{\halfsquare(0.4cm,0.2cm,#3)}}}

\def\spin4box(#1,#2,#3,#4){%
	\normalbaselines\m@th\offinterlineskip
	\vcenter{\hbox{\halfsquare(0.4cm,0.2cm,#1)}
	      \vskip-0.4pt
	      \hbox{\halfsquare(0.4cm,0.2cm,#2)}	
	      \vskip-0.4pt
              \hbox{\halfsquare(0.4cm,0.2cm,#3)}	
	      \vskip-0.4pt
	      \hbox{\halfsquare(0.4cm,0.2cm,#4)}}}

\def\minispinbox(#1,#2){%
	\normalbaselines\m@th\offinterlineskip
	\vcenter{\hbox{\halfsquare(0.4cm,0.3cm,#1)}
	      \vskip-0.4pt
	      \hbox{\halfsquare(1.0cm,0.3cm,\vdots)}	
	      \vskip-0.4pt
	      \hbox{\halfsquare(0.4cm,0.3cm,#2)}}}

\def\longbox(#1,#2){%
	\normalbaselines\m@th\offinterlineskip
	\vcenter{\hbox{\halfsquare(0.4cm,0.4cm,#1)}
	      \vskip-0.4pt
	      \hbox{\halfsquare(0.7cm,0.4cm,\vdots)}	
	      \vskip-0.4pt
	      \hbox{\halfsquare(0.4cm,0.4cm,#2)}}}

\def\threeseudobox(#1,#2){%
	\normalbaselines\m@th\offinterlineskip
	\vtop{\hbox{\topseudosquare(0.4cm,0.4cm,#1)}
	      \vskip-0.4pt
	      \hbox{\topseudosquare(0.8cm,0.4cm,\vdots)}	
	      \vskip-0.4pt
	      \hbox{\topseudosquare(0.4cm,0.4cm,#2)}}}

\def\topLongbox(#1,#2){%
	\normalbaselines\m@th\offinterlineskip
	\vtop{\hbox{\tophalfsquare(0.4cm,0.4cm,#1)}
	      \vskip-0.4pt
	      \hbox{\tophalfsquare(1.4cm,0.4cm,\vdots)}	
	      \vskip-0.4pt
	      \hbox{\tophalfsquare(0.4cm,0.4cm,#2)}}}

\def\MNbox(#1,#2,#3,#4){
       \vcenter{\hbox{
           $\topLongbox(#1,#3)
           \hskip-4pt \threeseudobox(#2,#4)$}}}

\def\spinMbox(#1,#2,#3,#4){
       \vcenter{\hbox{
           $\topspinbox(#1,#3)
           \hskip-4pt \threeseudobox(#2,#4)$}}}

\def\hfourbox(#1,#2,#3,#4){%
	\fsquare(0.3cm,#1)\addsquare(0.3cm,#2)\addsquare(0.3cm,#3)\addsquare(0.3cm,#4)\black-box}

\def\Hfourbox(#1,#2,#3,#4){%
	\Fsquare(0.4cm,#1)\Addsquare(0.4cm,#2)\Addsquare(0.4cm,#3)\Addsquare(0.4cm,#4)}

\def\hthreebox(#1,#2,#3){%
	\fsquare(0.3cm,#1)\addsquare(0.3cm,#2)\addsquare(0.3cm,#3)}

\def\htwobox(#1,#2){%
	\fsquare(0.3cm,#1)\addsquare(0.3cm,#2)}

\def\halfonebox(#1,#2){%
	\halfsquare(0.4cm,0.2cm,#1)\Addsquare(0.4cm,{#2})}

\def\halftwo(#1,#2,#3){%
	\halfsquare(0.4cm,0.2cm,#1)\Addsquare(0.4cm,#2)\Addsquare(0.4cm,#3)}

\def\halffour(#1,#2,#3,#4,#5){%
	\halfsquare(0.4cm,0.2cm,#1)\Addsquare(0.4cm,#2)\Addsquare(0.4cm,#3)
        \Addsquare(0.4cm,#4)\Addsquare(0.4cm,#5)}

\def\halfthree(#1,#2,#3,#4){%
	\halfsquare(0.4cm,0.2cm,#1)\Addsquare(0.4cm,#2)\Addsquare(0.4cm,#3)%
        \Addsquare(0.4cm,#4)}

\def\halfthreehalftwo(#1,#2,#3,#4,#5,#6,#7){%
	\normalbaselines\m@th\offinterlineskip
	\vcenter{\hbox{\halfthree({#1},{#2},{#3},{#4})}
              \vskip-0.4pt
	      \hbox{\halftwo({#5},{#6},{#7})}}}

\def\halfthreehalfone(#1,#2,#3,#4,#5,#6){%
	\normalbaselines\m@th\offinterlineskip
	\vcenter{\hbox{\halfthree({#1},{#2},{#3},{#4})}
              \vskip-0.4pt
	      \hbox{\halfonebox({#5},{#6})}}}

\def\halftwohalfone(#1,#2,#3,#4,#5){%
	\normalbaselines\m@th\offinterlineskip
	\vcenter{\hbox{\halftwo({#1},{#2},{#3})}
              \vskip-0.4pt
	      \hbox{\halfonebox({#4},{#5})}}}

\def\halftwohalf(#1,#2,#3,#4){%
	\normalbaselines\m@th\offinterlineskip
	\vcenter{\hbox{\halftwo({#1},{#2},{#3})}
              \vskip-0.4pt
	      \hbox{\halfsquare(0.4cm,0.2cm,{#4})}}}

\def\halftwohalftwo(#1,#2,#3,#4,#5,#6){%
	\normalbaselines\m@th\offinterlineskip
	\vcenter{\hbox{\halftwo({#1},{#2},{#3})}
              \vskip-0.4pt
	      \hbox{\halftwo({#4},{#5},{#6})}}}

\def\halftwohalfone(#1,#2,#3,#4,#5){%
	\normalbaselines\m@th\offinterlineskip
	\vcenter{\hbox{\halftwo({#1},{#2},{#3})}
              \vskip-0.4pt
	      \hbox{\halfonebox({#4},{#5})}}}

\def\halfonehalfone(#1,#2,#3,#4){%
	\normalbaselines\m@th\offinterlineskip
	\vcenter{\hbox{\halfonebox({#1},{#2})}
              \vskip-0.4pt
	      \hbox{\halfonebox({#3},{#4})}}}

\def\halftwohalfonehalfonebhalf(#1,#2,#3,#4,#5,#6,#7){%
	\normalbaselines\m@th\offinterlineskip
	\vcenter{\hbox{\halftwo(#1,#2,#3)}
	      \vskip-0.4pt
	      \hbox{\halfonebox(#4,#5)}	
	      \vskip-0.4pt
	      \hbox{\halfonebox(#6,#7)}	
	      \vskip-0.4pt
	      \hbox{\Blackhalfbox}}}

\def\halftwohalfhalfbhalf(#1,#2,#3,#4,#5){%
	\normalbaselines\m@th\offinterlineskip
	\vcenter{\hbox{\halftwo(#1,#2,#3)}
	      \vskip-0.4pt
	      \hbox{\halfsquare(0.4cm,0.2cm,#4)}	
	      \vskip-0.4pt
	      \hbox{\halfsquare(0.4cm,0.2cm,#5)}	
	      \vskip-0.4pt
	      \hbox{\Blackhalfbox}}}

\def\halfonehalfonehalfonebhalfone(#1,#2,#3,#4,#5,#6){%
	\normalbaselines\m@th\offinterlineskip
	\vcenter{\hbox{\halfonebox(#1,#2)}
	      \vskip-0.4pt
	      \hbox{\halfonebox(#3,#4)}	
	      \vskip-0.4pt
	      \hbox{\halfonebox(#5,#6)}	
	      \vskip-0.4pt
	      \hbox{\halfoneblackbox}}}

\def\halfonehalfonehalfbhalf(#1,#2,#3,#4,#5){%
	\normalbaselines\m@th\offinterlineskip
	\vcenter{\hbox{\halfonebox(#1,#2)}
	      \vskip-0.4pt
	      \hbox{\halfonebox(#3,#4)}	
	      \vskip-0.4pt
	      \hbox{\halfsquare(0.4cm,0.2cm,#5)}	
	      \vskip-0.4pt
	      \hbox{\Blackhalfbox}}}

\def\halfthreehalftwohalfone(#1,#2,#3,#4,#5,#6,#7,#8,#9){%
	\normalbaselines\m@th\offinterlineskip
	\vcenter{\hbox{\halfthree(#1,#2,#3,#4)}
	      \vskip-0.4pt
	      \hbox{\halftwo(#5,#6,#7)}	
	      \vskip-0.4pt
	      \hbox{\halfonebox(#8,#9)}}}

\def\half3half2halfhalf(#1,#2,#3,#4,#5,#6,#7,#8,#9){%
	\normalbaselines\m@th\offinterlineskip
	\vcenter{\hbox{\halfthree(#1,#2,#3,#4)}
	      \vskip-0.4pt
	      \hbox{\halftwo(#5,#6,#7)}	
	      \vskip-0.4pt
	      \hbox{\halfsquare(0.4cm,0.2cm,#8)}
              \vskip-0.4pt
              \hbox{\halfsquare(0.4cm,0.2cm,#9)}}}

\def\vfourbox(#1,#2,#3,#4){%
	\normalbaselines\m@th\offinterlineskip
	\vcenter{\hbox{\fsquare(0.3cm,#1)}
	      \vskip-0.4pt
	      \hbox{\fsquare(0.3cm,#2)}	
	      \vskip-0.4pt
	      \hbox{\fsquare(0.3cm,#3)}	
	      \vskip-0.4pt
	      \hbox{\fsquare(0.3cm,#4)}}}

\def\vfourhalfbox(#1,#2,#3,#4){%
	\normalbaselines\m@th\offinterlineskip
	\vcenter{\hbox{\halfsquare(0.4cm,0.2cm,#1)}
	      \vskip-0.4pt
	      \hbox{\halfsquare(0.4cm,0.2cm,#2)}	
	      \vskip-0.4pt
	      \hbox{\halfsquare(0.4cm,0.2cm,#3)}	
	      \vskip-0.4pt
	      \hbox{\halfsquare(0.4cm,0.2cm,#4)}}}

\def\twotwotwobbox(#1,#2,#3,#4,#5,#6){%
	\normalbaselines\m@th\offinterlineskip
	\vcenter{\hbox{\htwobox(#1,#2)}
	      \vskip-0.4pt
	      \hbox{\htwobox(#3,#4)}	
	      \vskip-0.4pt
	      \hbox{\htwobox(#5,#6)}	
	      \vskip-0.4pt
	      \hbox{\blackbox}}}

\def\Twotwotwobbox(#1,#2,#3,#4,#5,#6){%
	\normalbaselines\m@th\offinterlineskip
	\vcenter{\hbox{\Htwobox(#1,#2)}
	      \vskip-0.4pt
	      \hbox{\Htwobox(#3,#4)}	
	      \vskip-0.4pt
	      \hbox{\Htwobox(#5,#6)}	
	      \vskip-0.4pt
	      \hbox{\Blackbox}}}

\def\twotwotwotwobbox(#1,#2,#3,#4,#5,#6){%
	\normalbaselines\m@th\offinterlineskip
	\vcenter{\hbox{\htwobox(#1,#2)}
	      \vskip-0.4pt
	      \hbox{\htwobox(#3,#4)}	
	      \vskip-0.4pt
	      \hbox{\htwobox(#5,#6)}	
	      \vskip-0.4pt
	      \hbox{\twoblackbox}}}

\def\threetwotwobbox(#1,#2,#3,#4,#5,#6,#7){%
	\normalbaselines\m@th\offinterlineskip
	\vcenter{\hbox{\hthreebox(#1,#2,#3)}
	      \vskip-0.4pt
	      \hbox{\htwobox(#4,#5)}	
	      \vskip-0.4pt
	      \hbox{\htwobox(#6,#7)}	
	      \vskip-0.4pt
	      \hbox{\blackbox}}}

\def\threetwotwobox(#1,#2,#3,#4,#5,#6,#7){%
	\normalbaselines\m@th\offinterlineskip
	\vcenter{\hbox{\Hthreebox(#1,#2,#3)}
	      \vskip-0.4pt
	      \hbox{\Htwobox(#4,#5)}	
	      \vskip-0.4pt
	      \hbox{\Htwobox(#6,#7)}}}

\def\fourthreebox(#1,#2,#3,#4,#5,#6,#7){%
	\normalbaselines\m@th\offinterlineskip
	\vcenter{\hbox{\Hfourbox(#1,#2,#3,#4)}
	      \vskip-0.4pt
	      \hbox{\Hthreebox(#5,#6,#7)}}}

\def\threetwoonebox(#1,#2,#3,#4,#5,#6){%
	\normalbaselines\m@th\offinterlineskip
	\vcenter{\hbox{\Hthreebox(#1,#2,#3)}
	      \vskip-0.4pt
	      \hbox{\Htwobox(#4,#5)}	
	      \vskip-0.4pt
	      \hbox{\fsquare(0.4cm,#6)}}}

\def\threethreeonebox(#1,#2,#3,#4,#5,#6,#7){%
	\normalbaselines\m@th\offinterlineskip
	\vcenter{\hbox{\Hthreebox(#1,#2,#3)}
	      \vskip-0.4pt
	      \hbox{\Hthreebox(#4,#5,#6)}	
	      \vskip-0.4pt
	      \hbox{\fsquare(0.4cm,#7)}}}

\def\fourtwoonebox(#1,#2,#3,#4,#5,#6,#7){%
	\normalbaselines\m@th\offinterlineskip
	\vcenter{\hbox{\Hfourbox(#1,#2,#3,#4)}
	      \vskip-0.4pt
	      \hbox{\Htwobox(#5,#6)}	
	      \vskip-0.4pt
	      \hbox{\fsquare(0.4cm,#7)}}}

\def\twotwoonebbox(#1,#2,#3,#4,#5){%
	\normalbaselines\m@th\offinterlineskip
	\vcenter{\hbox{\htwobox(#1,#2)}
	      \vskip-0.4pt
	      \hbox{\htwobox(#3,#4)}	
	      \vskip-0.4pt
	      \hbox{\fsquare(0.3cm,#5)}
	      \vskip-0.4pt
	      \hbox{\blackbox}}}

\def\twooneonebbox(#1,#2,#3,#4){%
	\normalbaselines\m@th\offinterlineskip
	\vcenter{\hbox{\Htwobox(#1,#2)}
	      \vskip-0.4pt
	      \hbox{\fsquare(0.4cm,#3)}	
	      \vskip-0.4pt
	      \hbox{\fsquare(0.4cm,#4)}
	      \vskip-0.4pt
	      \hbox{\Blackbox}}}

\def\Vfourbox(#1,#2,#3,#4){%
	\normalbaselines\m@th\offinterlineskip
	\vcenter{\hbox{\Fsquare(0.4cm,#1)}
	      \vskip-0.4pt
	      \hbox{\Fsquare(0.4cm,#2)}	
	      \vskip-0.4pt
	      \hbox{\Fsquare(0.4cm,#3)}	
	      \vskip-0.4pt
	      \hbox{\Fsquare(0.4cm,#4)}}}

\def\vthreebox(#1,#2,#3){%
	\normalbaselines\m@th\offinterlineskip
	\vcenter{\hbox{\fsquare(0.3cm,#1)}
	      \vskip-0.4pt
	      \hbox{\fsquare(0.3cm,#2)}	
	      \vskip-0.4pt
	      \hbox{\fsquare(0.3cm,#3)}}}

\def\vthreehalfbox(#1,#2,#3){%
	\normalbaselines\m@th\offinterlineskip
	\vcenter{\hbox{\halfsquare(0.4cm,0.2cm,#1)}
	      \vskip-0.4pt
	      \hbox{\halfsquare(0.4cm,0.2cm,#2)}	
	      \vskip-0.4pt
	      \hbox{\halfsquare(0.4cm,0.2cm,#3)}}}

\def\vtwobox(#1,#2){%
	\normalbaselines\m@th\offinterlineskip
	\vcenter{\hbox{\fsquare(0.3cm,#1)}
	      \vskip-0.4pt
	      \hbox{\fsquare(0.3cm,#2)}}}

\def\Hthreebox(#1,#2,#3){%
	\Fsquare(0.4cm,#1)\Addsquare(0.4cm,#2)\Addsquare(0.4cm,#3)}

\def\Htwobox(#1,#2){%
	\Fsquare(0.4cm,#1)\Addsquare(0.4cm,#2)}

\def\Vthreebox(#1,#2,#3){%
	\normalbaselines\m@th\offinterlineskip
	\vcenter{\hbox{\Fsquare(0.4cm,#1)}
	      \vskip-0.4pt
	      \hbox{\Fsquare(0.4cm,#2)}	
	      \vskip-0.4pt
	      \hbox{\Fsquare(0.4cm,#3)}}}

\def\Vtwobox(#1,#2){%
	\normalbaselines\m@th\offinterlineskip
	\vcenter{\hbox{\Fsquare(0.4cm,#1)}
	      \vskip-0.4pt
	      \hbox{\Fsquare(0.4cm,#2)}}}

\def\twoone(#1,#2,#3){%
	\normalbaselines\m@th\offinterlineskip
	\vcenter{\hbox{\htwobox({#1},{#2})}
	      \vskip-0.4pt
	      \hbox{\fsquare(0.3cm,#3)}}}

\def\Twoone(#1,#2,#3){%
	\normalbaselines\m@th\offinterlineskip
	\vcenter{\hbox{\Htwobox({#1},{#2})}
	      \vskip-0.4pt
	      \hbox{\Fsquare(0.4cm,#3)}}}

\def\threeone(#1,#2,#3,#4){%
	\normalbaselines\m@th\offinterlineskip
	\vcenter{\hbox{\hthreebox({#1},{#2},{#3})}
	      \vskip-0.4pt
	      \hbox{\fsquare(0.3cm,#4)}}}

\def\Threeone(#1,#2,#3,#4){%
	\normalbaselines\m@th\offinterlineskip
	\vcenter{\hbox{\Hthreebox({#1},{#2},{#3})}
	      \vskip-0.4pt
	      \hbox{\Fsquare(0.4cm,#4)}}}

\def\Threetwo(#1,#2,#3,#4,#5){%
	\normalbaselines\m@th\offinterlineskip
	\vcenter{\hbox{\Hthreebox({#1},{#2},{#3})}
	      \vskip-0.4pt
	      \hbox{\Htwobox({#4},{#5})}}}

\def\threetwo(#1,#2,#3,#4,#5){%
	\normalbaselines\m@th\offinterlineskip
	\vcenter{\hbox{\hthreebox({#1},{#2},{#3})}
	      \vskip-0.4pt
	      \hbox{\htwobox({#4},{#5})}}}

\def\twotwo(#1,#2,#3,#4){%
	\normalbaselines\m@th\offinterlineskip
	\vcenter{\hbox{\htwobox({#1},{#2})}
	      \vskip-0.4pt
	      \hbox{\htwobox({#3},{#4})}}}

\def\Twotwo(#1,#2,#3,#4){%
	\normalbaselines\m@th\offinterlineskip
	\vcenter{\hbox{\Htwobox({#1},{#2})}
	      \vskip-0.4pt
	      \hbox{\Htwobox({#3},{#4})}}}

\def\twooneone(#1,#2,#3,#4){%
	\normalbaselines\m@th\offinterlineskip
	\vcenter{\hbox{\htwobox({#1},{#2})}
	      \vskip-0.4pt
	      \hbox{\fsquare(0.3cm,#3)}
	      \vskip-0.4pt
	      \hbox{\fsquare(0.3cm,#4)}}}

\def\Twooneone(#1,#2,#3,#4){%
	\normalbaselines\m@th\offinterlineskip
	\vcenter{\hbox{\Htwobox({#1},{#2})}
	      \vskip-0.4pt
	      \hbox{\Fsquare(0.4cm,#3)}
	      \vskip-0.4pt
	      \hbox{\Fsquare(0.4cm,#4)}}}

\def\Twooneonebone(#1,#2,#3,#4){%
	\normalbaselines\m@th\offinterlineskip
	\vcenter{\hbox{\Htwobox({#1},{#2})}
	      \vskip-0.4pt
	      \hbox{\Fsquare(0.4cm,#3)}
	      \vskip-0.4pt
	      \hbox{\Fsquare(0.4cm,#4)}
              \vskip-0.4pt
	      \hbox{\Blackbox}}}

\def\Twotwoone(#1,#2,#3,#4,#5){%
	\normalbaselines\m@th\offinterlineskip
	\vcenter{\hbox{\Htwobox({#1},{#2})}
	      \vskip-0.4pt
	      \hbox{\Htwobox({#3},{#4})}
              \vskip-0.4pt
	      \hbox{\Fsquare(0.4cm,{#5})}}}

\def\twotwoone(#1,#2,#3,#4,#5){%
	\normalbaselines\m@th\offinterlineskip
	\vcenter{\hbox{\htwobox({#1},{#2})}
	      \vskip-0.4pt
	      \hbox{\htwobox({#3},{#4})}
              \vskip-0.4pt
	      \hbox{\fsquare(0.3cm,{#5})}}}

\def\onehalfonehalfhalf(#1,#2,#3,#4,#5){%
	\normalbaselines\m@th\offinterlineskip
	\vcenter{\hbox{\halfonebox({#1},{#2})}
	      \vskip-0.4pt
	      \hbox{\halfonebox({#3},{#4})}
              \vskip-0.4pt
	      \hbox{\halfsquare(0.4cm,0.2cm,{#5})}}}

\def\halfonehalf(#1,#2,#3){%
	\normalbaselines\m@th\offinterlineskip
	\vcenter{\hbox{\halfonebox({#1},{#2})}
              \vskip-0.4pt
	      \hbox{\halfsquare(0.4cm,0.2cm,{#3})}}}

\def\Threeoneone(#1,#2,#3,#4,#5){%
	\normalbaselines\m@th\offinterlineskip
	\vcenter{\hbox{\Hthreebox({#1},{#2},{#3})}
	      \vskip-0.4pt
	      \hbox{\Fsquare(0.4cm,{#4})}
              \vskip-0.4pt
	      \hbox{\Fsquare(0.4cm,{#5})}}}

\def\threeoneone(#1,#2,#3,#4,#5){%
	\normalbaselines\m@th\offinterlineskip
	\vcenter{\hbox{\hthreebox({#1},{#2},{#3})}
	      \vskip-0.4pt
	      \hbox{\fsquare(0.3cm,#4)}
              \vskip-0.4pt
	      \hbox{\fsquare(0.3cm,#5)}}}

\def\a{\fsquare(0.3cm){1}\addsquare(0.3cm)(2)\addsquare(0.3cm)(3)}

\def\b{\hbox{%
	\normalbaselines\m@th\offinterlineskip
	\vcenter{\hbox{\fsquare(0.3cm){2}}\vskip-0.4pt\hbox{\fsquare(0.3cm){2}}}}}

\def\c{\hbox{\normalbaselines\m@th\offinterlineskip%
	\vcenter{\hbox{\a}\vskip-0.4pt\hbox{\b}}}}


\dimen1=0.4cm\advance\dimen1 by -0.8pt
\def\ffsquare#1{%
	\fsquare(0.4cm,\hbox{#1})}

\def\naga{%
	\hbox{$\vcenter to 0.4cm{\normalbaselines\m@th
	\hrule\vfil\hbox to 1.2cm{\hfill$\cdots$\hfill}\vfil\hrule}$}}

\def\vnaga{\normalbaselines\m@th\baselineskip0pt\offinterlineskip%
	\vrule\vbox to 1.2cm{\vskip7pt\hbox to \dimen1{$\hfil\vdots\hfil$}\vfil}\vrule}

\def\dhbox{\hbox{$\ffsquare 1 \naga \ffsquare N$}}

\def\dvbox{\hbox{\normalbaselines\m@th\baselineskip0pt\offinterlineskip\vbox{%
	  \hbox{$\ffsquare 1$}\vskip-0.4pt\hbox{$\vnaga$}\vskip-0.4pt\hbox{$\ffsquare N$}}}}


\def\addsquare(#1,#2){\hbox{$
	\dimen1=#1 \advance\dimen1 by -0.8pt
	\vcenter to #1{\hrule height0.4pt depth0.0pt%
	\hbox to #1{%
	\vbox to \dimen1{\vss%
	\hbox to \dimen1{\hss$\scriptstyle~#2~$\hss}%
	\vss}%
	\vrule width0.4pt}%
	\hrule height0.4pt depth0.0pt}$}}
\newcommand{\cI}{{\mathcal I}}
\newcommand{\cA}{{\mathcal A}}
\newcommand{\cB}{{\mathcal B}}
\newcommand{\cC}{{\mathcal C}}
\newcommand{\cD}{{\mathcal D}}
\newcommand{\cF}{{\mathcal F}}
\newcommand{\cH}{{\mathcal H}}
\newcommand{\cK}{{\mathcal K}}
\newcommand{\cL}{{\mathcal L}}
\newcommand{\cM}{{\mathcal M}}
\newcommand{\cN}{{\mathcal N}}
\newcommand{\cO}{{\mathcal O}}
\newcommand{\cS}{{\mathcal S}}
\newcommand{\cV}{{\mathcal V}}
\newcommand{\fra}{\mathfrak a}
\newcommand{\frb}{\mathfrak b}
\newcommand{\frc}{\mathfrak c}
\newcommand{\frd}{\mathfrak d}
\newcommand{\fre}{\mathfrak e}
\newcommand{\frf}{\mathfrak f}
\newcommand{\frg}{\mathfrak g}
\newcommand{\frh}{\mathfrak h}
\newcommand{\fri}{\mathfrak i}
\newcommand{\frj}{\mathfrak j}
\newcommand{\frk}{\mathfrak k}
\newcommand{\frI}{\mathfrak I}
\newcommand{\fm}{\mathfrak m}
\newcommand{\frn}{\mathfrak n}
\newcommand{\frp}{\mathfrak p}
\newcommand{\fq}{\mathfrak q}
\newcommand{\frr}{\mathfrak r}
\newcommand{\frs}{\mathfrak s}
\newcommand{\frt}{\mathfrak t}
\newcommand{\fru}{\mathfrak u}
\newcommand{\frA}{\mathfrak A}
\newcommand{\frB}{\mathfrak B}
\newcommand{\frF}{\mathfrak F}
\newcommand{\frG}{\mathfrak G}
\newcommand{\frH}{\mathfrak H}
\newcommand{\frJ}{\mathfrak J}
\newcommand{\frN}{\mathfrak N}
\newcommand{\frP}{\mathfrak P}
\newcommand{\frT}{\mathfrak T}
\newcommand{\frU}{\mathfrak U}
\newcommand{\frV}{\mathfrak V}
\newcommand{\frX}{\mathfrak X}
\newcommand{\frY}{\mathfrak Y}
\newcommand{\frZ}{\mathfrak Z}
\newcommand{\rA}{\mathrm{A}}
\newcommand{\rC}{\mathrm{C}}
\newcommand{\rd}{\mathrm{d}}
\newcommand{\rB}{\mathrm{B}}
\newcommand{\rD}{\mathrm{D}}
\newcommand{\rE}{\mathrm{E}}
\newcommand{\rH}{\mathrm{H}}
\newcommand{\rK}{\mathrm{K}}
\newcommand{\rL}{\mathrm{L}}
\newcommand{\rM}{\mathrm{M}}
\newcommand{\rN}{\mathrm{N}}
\newcommand{\rR}{\mathrm{R}}
\newcommand{\rT}{\mathrm{T}}
\newcommand{\rZ}{\mathrm{Z}}
\newcommand{\bbA}{\mathbb A}
\newcommand{\bbC}{\mathbb C}
\newcommand{\bbG}{\mathbb G}
\newcommand{\bbF}{\mathbb F}
\newcommand{\bbH}{\mathbb H}
\newcommand{\bbP}{\mathbb P}
\newcommand{\bbN}{\mathbb N}
\newcommand{\bbQ}{\mathbb Q}
\newcommand{\bbR}{\mathbb R}
\newcommand{\bbV}{\mathbb V}
\newcommand{\bbZ}{\mathbb Z}
\newcommand{\adj}{\operatorname{adj}}
\newcommand{\Ad}{\mathrm{Ad}}
\newcommand{\Ann}{\mathrm{Ann}}
\newcommand{\rcris}{\mathrm{cris}}
\newcommand{\ch}{\mathrm{ch}}
\newcommand{\coker}{\mathrm{coker}}
\newcommand{\diag}{\mathrm{diag}}
\newcommand{\Diff}{\mathrm{Diff}}
\newcommand{\Dist}{\mathrm{Dist}}
\newcommand{\rDR}{\mathrm{DR}}
\newcommand{\ev}{\mathrm{ev}}
\newcommand{\Ext}{\mathrm{Ext}}
\newcommand{\cExt}{\mathcal{E}xt}
\newcommand{\fin}{\mathrm{fin}}
\newcommand{\Frac}{\mathrm{Frac}}
\newcommand{\GL}{\mathrm{GL}}
\newcommand{\Hom}{\mathrm{Hom}}
\newcommand{\hd}{\mathrm{hd}}
\newcommand{\rht}{\mathrm{ht}}
\newcommand{\id}{\mathrm{id}}
\newcommand{\im}{\mathrm{im}}
\newcommand{\inc}{\mathrm{inc}}
\newcommand{\ind}{\mathrm{ind}}
\newcommand{\coind}{\mathrm{coind}}
\newcommand{\Lie}{\mathrm{Lie}}
\newcommand{\Max}{\mathrm{Max}}
\newcommand{\mult}{\mathrm{mult}}
\newcommand{\op}{\mathrm{op}}
\newcommand{\ord}{\mathrm{ord}}
\newcommand{\pt}{\mathrm{pt}}
\newcommand{\qt}{\mathrm{qt}}
\newcommand{\rad}{\mathrm{rad}}
\newcommand{\res}{\mathrm{res}}
\newcommand{\rgt}{\mathrm{rgt}}
\newcommand{\rk}{\mathrm{rk}}
\newcommand{\SL}{\mathrm{SL}}
\newcommand{\soc}{\mathrm{soc}}
\newcommand{\Spec}{\mathrm{Spec}}
\newcommand{\St}{\mathrm{St}}
\newcommand{\supp}{\mathrm{supp}}
\newcommand{\Tor}{\mathrm{Tor}}
\newcommand{\Tr}{\mathrm{Tr}}
\newcommand{\wt}{\mathrm{wt}}
\newcommand{\Ab}{\mathbf{Ab}}
\newcommand{\Alg}{\mathbf{Alg}}
\newcommand{\Grp}{\mathbf{Grp}}
\newcommand{\Mod}{\mathbf{Mod}}
\newcommand{\Sch}{\mathbf{Sch}}\newcommand{\bfmod}{{\bf mod}}
\newcommand{\Qc}{\mathbf{Qc}}
\newcommand{\Rng}{\mathbf{Rng}}
\newcommand{\Top}{\mathbf{Top}}
\newcommand{\Var}{\mathbf{Var}}
\newcommand{\gromega}{\langle\omega\rangle}
\newcommand{\lbr}{\begin{bmatrix}}
\newcommand{\rbr}{\end{bmatrix}}
\newcommand{\for}{\bigcirc\kern-2.6ex \because}
\newcommand{\forb}{\bigcirc\kern-2.8ex \because}
\newcommand{\forbb}{\bigcirc\kern-3.0ex \because}
\newcommand{\forbbb}{\bigcirc\kern-3.1ex \because}
\newcommand{\cd}{commutative diagram }
\newcommand{\SpS}{spectral sequence}
\newcommand\C{\mathbb C}
\newcommand\hh{{\hat{H}}}
\newcommand\eh{{\hat{E}}}
\newcommand\F{\mathbb F}
\newcommand\fh{{\hat{F}}}

\def\AA{{\cal A}}
\def\al{\alpha}
\def\bq{B_q(\ge)}
\def\bqm{B_q^-(\ge)}
\def\bqz{B_q^0(\ge)}
\def\bqp{B_q^+(\ge)}
\def\beneme{\begin{enumerate}}
\def\beq{\begin{equation}}
\def\beqn{\begin{eqnarray}}
\def\beqnn{\begin{eqnarray*}}
\def\bigsl{{\hbox{\fontD \char'54}}}
\def\bbra#1,#2,#3{\left\{\begin{array}{c}\hspace{-5pt}
#1;#2\\ \hspace{-5pt}#3\end{array}\hspace{-5pt}\right\}}
\def\cd{\cdots}
\def\CC{\hbox{\bf C}}
\def\ddd{\hbox{\germ D}}
\def\del{\delta}
\def\Del{\Delta}
\def\Delr{\Delta^{(r)}}
\def\Dell{\Delta^{(l)}}
\def\Delb{\Delta^{(b)}}
\def\Deli{\Delta^{(i)}}
\def\Delre{\Delta^{\rm re}}
\def\ei{e_i}
\def\eit{\tilde{e}_i}
\def\eneme{\end{enumerate}}
\def\ep{\epsilon}
\def\eeq{\end{equation}}
\def\eeqn{\end{eqnarray}}
\def\eeqnn{\end{eqnarray*}}
\def\fit{\tilde{f}_i}
\def\FF{{\rm F}}
\def\ft{\tilde{f}}
\def\gau#1,#2{\left[\begin{array}{c}\hspace{-5pt}#1\\
\hspace{-5pt}#2\end{array}\hspace{-5pt}\right]}
\def\ge{\hbox{\germ g}}
\def\gl{\hbox{\germ gl}}
\def\hom{{\hbox{Hom}}}
\def\ify{\infty}
\def\io{\iota}
\def\kp{k^{(+)}}
\def\km{k^{(-)}}
\def\llra{\relbar\joinrel\relbar\joinrel\relbar\joinrel\rightarrow}
\def\lan{\langle}
\def\lar{\longrightarrow}
\def\max{{\rm max}}
\def\lm{\lambda}
\def\Lm{\Lambda}
\def\mapright#1{\smash{\mathop{\longrightarrow}\limits^{#1}}}
\def\mm{{\bf{\rm m}}}
\def\nd{\noindent}
\def\nn{\nonumber}
\def\nnn{\hbox{\germ n}}
\def\catob{{\cal O}(B)}
\def\oint{{\cal O}_{\rm int}(\ge)}
\def\ot{\otimes}
\def\op{\oplus}
\def\opi{\ovl\pi_{\lm}}
\def\ovl{\overline}
\def\plm{\Psi^{(\lm)}_{\io}}
\def\qq{\qquad}
\def\q{\quad}
\def\qed{\hfill\framebox[2mm]{}}
\def\QQ{\hbox{\bf Q}}
\def\qi{q_i}
\def\qii{q_i^{-1}}
\def\ra{\rightarrow}
\def\ran{\rangle}
\def\rlm{r_{\lm}}
\def\ssl{\hbox{\germ sl}}
\def\slh{\widehat{\ssl_2}}
\def\ti{t_i}
\def\tii{t_i^{-1}}
\def\til{\tilde}
\def\tm{\times}
\def\tt{{\hbox{\germ{t}}}}
\def\ttt{\hbox{\germ t}}
\def\ua{U_{\AA}}
\def\ue{U_{\vep}}
\def\uq{U_q(\ge)}
\def\ufin{U^{\rm fin}_{\vep}}
\def\ufinp{(U^{\rm fin}_{\vep})^+}
\def\ufinm{(U^{\rm fin}_{\vep})^-}
\def\ufinz{(U^{\rm fin}_{\vep})^0}
\def\uqm{U^-_q(\ge)}
\def\uqp{U^+_q(\ge)}
\def\uqmq{{U^-_q(\ge)}_{\bf Q}}
\def\uqpm{U^{\pm}_q(\ge)}
\def\uqq{U_{\bf Q}^-(\ge)}
\def\uqz{U^-_{\bf Z}(\ge)}
\def\ures{U^{\rm res}_{\AA}}
\def\urese{U^{\rm res}_{\vep}}
\def\uresez{U^{\rm res}_{\vep,\ZZ}}
\def\util{\widetilde\uq}
\def\uup{U^{\geq}}
\def\ulow{U^{\leq}}
\def\bup{B^{\geq}}
\def\blow{\ovl B^{\leq}}
\def\vep{\varepsilon}
\def\vp{\varphi}
\def\vpi{\varphi^{-1}}
\def\VV{{\cal V}}
\def\xii{\xi^{(i)}}
\def\Xiioi{\Xi_{\io}^{(i)}}
\def\WW{{\cal W}}
\def\wtil{\widetilde}
\def\what{\widehat}
\def\wpi{\widehat\pi_{\lm}}
\def\ZZ{\mathbb Z}

\renewcommand{\thesection}{\arabic{section}}
\section{Introduction}
\setcounter{equation}{0}
\renewcommand{\theequation}{\thesection.\arabic{equation}}

The notion of crystals is initiated by Kashiwara 
(\cite{K1},\cite{K3},\cite{K5}), which influences over many areas
in mathematics, in particular, 
combinatorics and 
 representation theory, {\it e.g.,} combinatorics
of Young tableaux($=$ semi-standard tableaux), 
piece-wise linear combinatorics, {\it etc.}
Indeed, in \cite{KN}, we succeed in describing the crystal bases
for classical quantum algebras by using Young tableaux.
One feature of crystal theory
 is that it produces many piece-wise
linear formulae (\cite{K5},\cite{N1},\cite{N2},\cite{NZ}).

Theory of geometric crystals is introduced by Berenstein
and Kazhdan \cite{BK} in semi-simple setting and is extended 
to Kac-Moody setting in \cite{N}, which is a kind of 
geometric analogue of Kashiwara's crystal theory.
More precisely, let $G$ be a Kac-Moody group over $\bbC$, 
$T$ be its maximal torus and $I$ 
be a finite index set of its simple roots.
For an ind-(algebraic)variety $X$, 
morphisms $e_i:\bbC^\times\times X\rightarrow X$ 
$(i\in I)$ and $\gamma:X\rightarrow T$, the triplet
$(X, \gamma, \{e_i\}_{i\in I})$ is called 
a {\it geometric crystal}
if they satisfy the conditions as in 
Definition \ref{geo-cry}. 
Geometric crystals are not only analogy of 
crystals, but also has certain categorical 
correpondence to crystals, which is called a
tropicalization/ultra-discretization.
It is so remarkable that this 
correpondence reproduces several piece-wise linear
formulae in the theory of crystals from 
subtraction free(=positive) rational formulae in geometric crystals 
(\cite{N}) as follows:
\[
\xymatrix@R=3pt{
\{\rm Geometric\,\, Crystals\}\ar@{->}^{
\rm\scriptstyle \,\,\,\,\,\,\,\,ultra-discretization} 
@<1pt>[rr]
\ar@{<-}_{
\rm\scriptstyle \,\,\,\,tropicalization} @<-3pt>[rr]
&& \{\rm Crystals\} \\
x\times y,\,\,x/y, \,\,x+y &&
x+y,\,\,x-y,\,\,{\rm max}(x,y)
}
\]
Furthermore, this correpondence reproduces the tensor product 
structure of crystals from the product structure 
of geometric crystals (\cite{BK}).

Let $B$ be a Borel subgroup
of $G$ and $W$ be the Weyl group associated with $G$.
Any finite Schubert variety $\ovl X_w\subset X:=G/B$ has 
a natural geometric crystal structure(\cite{BK},\cite{N}).
Then, in semi-simple setting we know that 
the whole flag variety $X:=G/B$
holds a geometric crystal structure. But, 
in general Kac-Moody setting, we do not have any natural geometric
crystal structure on the flag variety $X$.
The opposite unipotent subgroup $U^-$ can be seen as an 
open dense subset of $X$.
In this paper, we present some sufficient condition 
for existence of geometric(unipotent) crystal structure
on $U^-$ and then on $X$, which is described as follows:
if there exists a morphism ${\cal T}:U^-\rightarrow T$
satisfying the condition as in Lemma \ref{suf}, then we obtain 
$U$-morphism $F:U^-\rightarrow B^-$ and then 
the associated unipotent crystal structure, which means the exsistence of 
a geometric crystal structure on $U^-$.
In semi-simple cases, there exists such morphism which is given by 
matrix coefficients. 
In particular, for $G=SL_{n+1}(\bbC)$ case, we present
its geometric crystal structure explicitly
and reveal that it corresponds to the crystals
called {\it generalized Young tableaux}, which is a 
sort of ``limit'' of usual Young tableaux and 
forms a free $\ZZ$-lattice of rank $\frac{n(n+1)}{2}$. 
In more general cases, {\it e.g.,} affine cases, 
the existence
of such morphisms is not yet known, which is our 
further problem.

The article is organized as follows:
in Sect.2, we review the notion of geometric crystals, 
unipotent crystals and the 
tropicalization/ultra-discretization correspondence.
In Sect.3, we consider geometric crystal on a
unipotent subgroup $U^-\subset G$ and in Sect.4, the explicit 
geometric crystal strcutre on $U^-\subset SL_{n+1}(\bbC)$
is described. In the final section, 
we give a tropcialization/ultra-discretization correpondence
between geometric crystals on $U^-$ and generalized Young 
tableaux.

The author would acknowledge M.Kashiwara for valuable 
discussions and comments.

\renewcommand{\thesection}{\arabic{section}}
\section{Geometric Crystals and Unipotent Crystals}
\setcounter{equation}{0}
\renewcommand{\theequation}{\thesection.\arabic{equation}}

\newtheorem{pro2}{Proposition}[section]
\theoremstyle{definition}
\newtheorem{def2}[pro2]{Definition}
\theoremstyle{plain}
\newtheorem{lem2}[pro2]{Lemma}
\newtheorem{thm2}[pro2]{Theorem}
\newtheorem{ex2}[pro2]{Example}

\subsection{Kac-Moody algebras and Kac-Moody groups}
Fix a symmetrizable generalized Cartan matrix
 $A=(a_{ij})_{i,j\in I}$, where $I$ be a finite index set.
Let $(\tt,\{\al_i\}_{i\in I},\{h_i\}_{i\in I})$ 
be the associated
root data, where ${\tt}$ be the vector space 
over $\bbC$ with 
dimension $|I|+$ corank$(A)$, and 
the set of simple roots $\{\al_i\}_{i\in I}\subset\tt^*$ and 
the set of simple co-roots $\{h_i\}_{i\in I}\subset\tt$
are linearly independent indexed sets 
satisfying $\al_i(h_j)=a_{ij}$.

The Kac-Moody Lie algebra $\ge=\ge(A)$ associated with $A$
is the Lie algebera over $\bbC$ generated by $\tt$, the 
Chevalley generators $e_i$ and $f_i$ $(i\in I)$
with the usual defining relations (\cite{KP},\cite{PK}).
There is the root space decomposition 
$\ge=\bigoplus_{\al\in \tt^*}\ge_{\al}$.
Denote the set of roots by 
$\Delta:=\{\al\in \tt^*|\al\ne0,\,\,\ge_{\al}\ne(0)\}$.
Set $Q:=\sum_i\bbZ \al_i$, $Q_+:=\sum_i\bbZ_{\geq0} \al_i$
and $\Delta_+:=\Delta\cap Q_+$.
An element of $\Delta_+$ is called a positive root.
Let $\omega$ be the Chevalley involution of $\ge$ defined by
$\omega(e_i)=-f_i$, $\omega(f_i)=-e_i$ and $\omega(h)=-h$
for $h\in \tt$.
Let $L(\Lm)$ $(\Lm\in P_+:\text{set of dominant weights})$ 
be an irreducible integrable highest weight module with 
highest weight $\Lm$ and 
$\pi_\Lm:\ge\rightarrow {\rm End}(L(\Lm))$ 
be the $\ge$-action. 
The action $\pi^*_\Lm:=\pi_\Lm\circ \omega$ defines a 
$\ge$-module structure on $L(\Lm)$, 
which is called the contragredient
module of $L(\Lm)$ and denoted $L^*(\Lm)$.
Let us fix a highest weight vector
$u_\Lm\in L(\lm)$ and denote it by $u^*_\Lm$ in $L^*(\Lm)$.
We obtain a unique 
$\ge$-invariant bilinear form 
$\lan\q,\q\ran$ on $L(\Lm)\times L^*(\Lm)$ such that 
$\lan u_\Lm,u^*_\Lm\ran=1$.

Define simple reflections $s_i\in{\rm Aut}(\tt)$ $(i\in I)$ by
$s_i(h):=h-\al_i(h)h_i$, which generate the Weyl group $W$.
We also define the action of $W$ on $\tt^*$ by
$s_i(\lm):=\lm-\al(h_i)\al_i$.
Set $\Delre:=\{w(\al_i)|w\in W,\,\,i\in I\}$, whose element 
is called a real root.

Let $G$ be the Kac-Moody group associated with the 
derived Lie algebra $\ge'$ defined in \cite{PK}.
Set  $U_{\al}:=\exp\ge_{\al}$ $(\al\in \Delre)$,
which is an one-parameter subgroup of $G$ and 
$G$ is generated by $U_{\al}$ $(\al\in \Delre)$.
Let $U^{\pm}$ be the subgroups generated by $U_{\pm\al}$
($\al\in \Delre_+=\Delre\cap Q_+$), {\it i.e.,}
$U^{\pm}:=\lan U_{\pm\al}|\al\in\Del^{\rm re}_+\ran$,
which is  called the unipotent subgroup of $G$.
Here note that if $\ge$ is a semi-simple Lie algebra,
then $G$ is a usual semi-simple algebraic group over $\bbC$.

For any $i\in I$, there exists a unique homomorphism;
$\phi_i:SL_2(\bbC)\rightarrow G$ such that
\[
x_i(t):=\phi_i\left(
\left(
\begin{array}{cc}
1&t\\
0&1
\end{array}
\right)\right)=\exp t e_i,\,\,
y_i(t):= \phi_i\left(
\left(
\begin{array}{cc}
1&0\\
t&1
\end{array}
\right)\right)=\exp t f_i\,\,\,\,(t\in\bbC).
\]
Set $G_i:=\phi_i(SL_2(\bbC))$,
$T_i:=\phi_i(\{{\rm diag}(t,t^{-1})|t\in\bbC\})$ and 
$N_i:=N_{G_i}(T_i)$. Let
$T$ (resp. $N$) be the subgroup of $G$ generated by $T_i$
(resp. $N_i$), which is called a {\it maximal torus} in $G$ and
$B^{\pm}:=U^{\pm}T$ be the Borel subgroup of $G$.
We have the isomorphism
$\phi:W\mapright{\sim}N/T$ defined by $\phi(s_i)=N_iT/T$.
An element $\ovl s_i:=x_i(-1)y_i(1)x_i(-1)$ is in 
$N_G(T)$, which is a representative of 
$s_i\in W=N_G(T)/T$. 
Define $R(w)$ for $w\in W$ by
\[
 R(w):=\{(i_1,i_2,\cd,i_l)\in I^l|w=s_{i_1}s_{i_2}\cd s_{i_l}\},
\]
where $l$ is the length of $w$.
We associate to each $w\in W$ its standard representative 
$\bar w\in N_G(T)$ by 
$\bar w=\bar s_{i_1}\bar s_{i_2}\cd\bar s_{i_l}$
for any $(i_1,i_2,\cd,i_l)\in R(w)$.

We have the following (as for ind-variety and ind-group, 
see \cite{Ku2}):
\begin{pro2}[\cite{Ku2}]
\label{ind-g}
\begin{enumerate}
\item
Let $G$ be a Kac-Moody group and 
$U^\pm$, $B^\pm$ be its subgroups as above.
Then $G$ is an ind-group and $U^\pm$, $B^\pm$
are its closed ind-subgroups.
\item
The multiplication maps 
\[
\begin{array}{ccc}
T\tm U &\longrightarrow &B \\
(t,u)& \mapsto & tu
\end{array}\qq
\begin{array}{ccc}
U^-\tm T&\longrightarrow & B^-\\
(v,t)&\mapsto& vt
\end{array}
\]
are isomorphisms of ind-varieties.
\end{enumerate}
\end{pro2}
\subsection{Geometric Crystals}
In this subsection, we review the notion of geometric crystals
(\cite{BK},\cite{N}).

Let $(a_{ij})_{i,j\in I}$ be a symmetrizable
generalized Cartan matrix and 
$G$ be the associated Kac-Moody group with the
maximal torus $T$.
An element in  $\hom(T,\bbC^\times)$
(resp. $\hom(\bbC^\times,T)$) 
is called a {\it character}
(resp. {\it co-character}) of $T$.
We define a {\it simple co-root} 
$\al_i^\vee\in \hom(\bbC^\times,T)$
$(i\in I)$ by $\al_i^\vee(t):=T_i$. 
We have a pairing $\lan \al^\vee_i,\al_j\ran=a_{ij}$.

Let $X$ be an ind-variety over $\bbC$, 
$\gamma:X\rightarrow T$
be a rational morphism and a family of 
rational morphisms 
$e_i:\bbC^\times \times X\rightarrow X$ 
$(i\in I)$; 
\[
\begin{array}{cccc}
 e^c_i&:\bbC^\times \times X&\longrightarrow& X\\
&(c,x)&\mapsto& e^c_i(x).
\end{array}
\]
For a word ${\bf i}=(i_1,i_2,\cd,i_l)\in R(w)$ 
$(w\in W)$,
set $\al^{(l)}:=\al_{i_l}$, 
$\al^{(l-1)}:=s_{i_l}(\al_{i_{l-1}})$, $\cd$,
$\al^{(1)}:=s_{i_l}\cd s_{i_2}(\al_{i_1})$.
Now for a word ${\bf i}=(i_1,i_2,\cd,i_l)\in R(w)$
we define a rational morphism $e_{\bf i}:T\times X\rightarrow X$
by
\[
(t,x)\mapsto e_{\bf i}^t(x):=e_{i_1}^{\al^{(1)}(t)}
e_{i_2}^{\al^{(2)}(t)}\cd e_{i_l}^{\al^{(l)}(t)}(x).
\]

\begin{def2}
\label{geo-cry}
\begin{enumerate}
\item
The triplet $\chi=(X,\gamma,\{e_i\}_{i\in I})$
is a {\it geometric crystal} if it 
satisfies $e^1(x)=x$ and 
\begin{eqnarray}
&&\gamma(e^c_i(x))=\al_i^\vee(c)\gamma(x),\label{gamma}\\
&&e_{\bf i}=e_{\bf i'}
\q\text{for any $w\in W$, and any ${\bf i}$, 
${\bf i'}\in R(w)$}.
\label{ei=ei'}
\end{eqnarray}
\item
Let $(X,\gamma_X,\{e^X_i\}_{i\in I})$ and 
$(Y,\gamma_Y,\{e^Y_i\}_{i\in I})$ be geometric crystals. 
A rational morphism $f:X\ra Y$ is a 
{\it morphism of geometric crystals}
if $f$ satisfies that 
\[
f\circ e^X_i=e^Y_i\circ f, \q 
\gamma_X=\gamma_Y\circ f.
\]
In particular, if a morphism $f$ is a birational isomorphism 
of ind-varieties, it is called an 
{\it isomorphism of geometric crystals}.
\end{enumerate}
\end{def2}
The following lemma is a direct result from 
\cite{BK}[Lemma 2.1] and the fact that 
the Weyl group of any 
Kac-Moody Lie algebra is a Coxeter group 
\cite{Kac}[Proposition 3.13].
\begin{lem2}
\label{Verma}
The relations $(\ref{ei=ei'})$ 
are equivalent to the following
relations:
\[
 \begin{array}{lll}
&\hspace{-20pt}e^{c_1}_{i}e^{c_2}_{j}
=e^{c_2}_{j}e^{c_1}_{i}&
{\rm if }\,\,\lan \al^\vee_i,\al_j\ran=0,\\
&\hspace{-20pt} e^{c_1}_{i}e^{c_1c_2}_{j}e^{c_2}_{i}
=e^{c_2}_{j}e^{c_1c_2}_{i}e^{c_1}_{j}&
{\rm if }\,\,\lan \al^\vee_i,\al_j\ran
=\lan \al^\vee_j,\al_i\ran=-1,\\
&\hspace{-20pt}
e^{c_1}_{i}e^{c^2_1c_2}_{j}e^{c_1c_2}_{i}e^{c_2}_{j}
=e^{c_2}_{j}e^{c_1c_2}_{i}e^{c^2_1c_2}_{j}e^{c_1}_{i}&
{\rm if }\,\,\lan \al^\vee_i,\al_j\ran=-2,\,
\lan \al^\vee_j,\al_i\ran=-1,\\
&\hspace{-20pt}
e^{c_1}_{i}e^{c^2_1c_2}_{j}e^{c^3_1c_2}_{i}
e^{c^3_1c^2_2}_{j}e^{c_1c_2}_{i}e^{c_2}_{j}
=e^{c_2}_{j}e^{c_1c_2}_{i}e^{c^3_1c^2_2}_{j}e^{c^3_1c_2}_{i}
e^{c^2_1c_2}_je^{c_1}_i&
{\rm if }\,\,\lan \al^\vee_i,\al_j\ran=-3,\,
\lan \al^\vee_j,\al_i\ran=-1,
\end{array}
\]
\end{lem2}
{\sl Remark.} 
If $\lan\al^\vee_i,\al_j\ran\lan \al^\vee_j,\al_i\ran
\geq4$, there is no relation between $e_i$ and $e_j$.

\subsection{Unipotent Crystals}
\label{uni-cry}

In the sequel, we denote the unipotent subgroup 
$U^+$ by $U$. 
We define unipotent crystals (see \cite{BK}) associated to 
Kac-Moody groups. 
The definitions below follow \cite{BK},\cite{N}.
\begin{def2}
Let $X$ be an ind-variety over $\bbC$ and 
$\al:U\times X\rightarrow X$ be a rational $U$-action
such that $\al$ is defined on $\{e\}\times X$. Then, 
the pair ${\bf X}=(X,\al)$ is called a $U$-{\it variety}. 
For $U$-varieties ${\bf X}=(X,\al_X)$
and ${\bf Y}=(Y,\al_Y)$, 
a rational morphism
$f:X\rightarrow Y$ is called a 
$U$-{\it morphism} if it commutes with
the action of $U$.
\end{def2}
Now, we define the $U$-variety structure on $B^-=U^-T$.
By Proposition \ref{ind-g},
 $B^-$ is an ind-subgroup of $G$ and then
is an ind-variety over $\bbC$.
The multiplication map in $G$ induces the open embedding;
$ B^-\times U\hookrightarrow G,$
then this is a birational isomorphism. 
Let us denote the inverse birational isomorphism by $g$;
\[
 g:G\longrightarrow B^-\times U.
\]
Then we define the rational morphisms 
$\pi^-:G\rightarrow B^-$ and 
$\pi:G\rightarrow U$ by 
$\pi^-:={\rm proj}_{B^-}\circ g$ 
and $\pi:={\rm proj}_U\circ g$.
Now we define the rational $U$-action $\al_{B^-}$ on $B^-$ by 
\[
 \al_{B^-}:=\pi^-\circ m:U\times B^-\longrightarrow B^-,
\]
where $m$ is the multiplication map in $G$.
Then we obtain $U$-variety ${\bf B}^-=(B^-,\al_{B^-})$.
\begin{def2}
\label{uni-def}
\begin{enumerate}
\item
Let ${\bf X}=(X,\al)$ 
be a $U$-variety and $f:X \rightarrow {\bf B^-}$ 
be a $U$-morphism.
The pair $({\bf X}, f)$ is called 
a {\it unipotent $G$-crystal}
or, for short, {\it unipotent crystal}.
\item
Let $({\bf X},f_X)$ and $({\bf Y},f_Y)$ 
be unipotent crystals.
A $U$-morphism $g:X\ra Y$ is called a {\it morphism of 
unipotent crystals} if $f_X=f_Y\circ g$.
In particular, if $g$ is a birational isomorphism
of ind-varieties, it is called an {\it isomorphism of 
unipotent crystals}.
\end{enumerate}
\end{def2}
We define a product of 
unipotent crystals following \cite{BK}.
For unipotent crystals $({\bf X},f_X)$, $({\bf Y},f_Y)$, 
define a morphism 
$\al_{X\times Y}:U\tm X\tm Y\rightarrow X\tm Y$ by
\begin{equation}
\al_{X\tm Y}(u,x,y):=(\al_X(u,x),\al_Y(\pi(u\cdot f_X(x)),y)).
\label{XY}
\end{equation}
If there is no confusion,
we use abbreviated notation $u(x,y)$ 
for $\al_{X\tm Y}(u,x,y)$.
\begin{thm2}[\cite{BK}]
\label{prod}
\begin{enumerate}
\item
The morphism $\al_{X\tm Y}$ defined above 
is a rational $U$-morphism
on $X\tm Y$.
\item
Let ${\bf m}:B^-\tm B^-\rightarrow B^-$ 
be a multiplication morphism 
and $f=f_{X\tm Y}:X\tm Y\rightarrow B^-$ be the 
rational morphism defined by 
\[
f_{X\tm Y}:={\bf m}\circ( f_X\tm f_Y).
\]
Then $ f_{X\tm Y}$ is a $U$-morphism and then, 
$({\bf X\tm Y}, f_{X\tm Y})$ is a unipotent crystal, 
which we call a product of unipotent crystals  
$({\bf X},f_X)$ and $({\bf Y},f_Y)$.
\item 
Product of unipotent crystals is associative.
\end{enumerate}
\end{thm2}
\subsection{From unipotent crystals to geometric crystals}
\label{u->g}
We have the canonical projection 
$\xi_i:U^-\rightarrow U_{-\al_i}$ $(i\in I)$ (see \cite{N}).
Now, we define the function on $U^-$ by 
\[
\chi_i:=y_i^{-1}\circ\xi_i:
U^-\longrightarrow U_{-\al_i}\longrightarrow
\bbC,
\]
and extend this to the function on 
$B^-$ by $\chi_i(u\cdot t):=\chi_i(u)$ for 
$u\in U^-$ and $t\in T$.
For a unipotent $G$-crystal $\bf(X,f_X)$, we define a function
$\vp_i:=\vp_i^X:X\rightarrow \bbC$ by 
\[
\vp_i:=\chi_i\circ{\bf f_X},
\]
and a rational morphism 
$\gamma_X:X\rightarrow T$ by 
\begin{equation}
\gamma_{X}:=
{\rm proj}_T\circ{\bf f_X}:X\rightarrow B^-\rightarrow T,
\label{gammax}
\end{equation}
where ${\rm proj}_T$ is the canonical projection.
Suppose that the function $\vp_i$ is not identically zero on $X$. 
We define a rational 
morphism $e_i:\bbC^\tm\tm X\rightarrow X$ by
\begin{equation}
e^c_i(x):=x_i
\left({\frac{c-1}{\vp_i(x)}}\right)(x).
\label{ei}
\end{equation}
\begin{thm2}[\cite{BK}]
\label{U-G}
For a unipotent $G$-crystal $\bf(X,f_X)$, 
suppose that 
the function $\vp_i$ is not identically zero
for any $i\in I$.
Then the rational morphisms $\gamma_X:X\rightarrow T$ 
and 
$e_i:\bbC^\tm\tm X\rightarrow X$ as above 
define a geometric 
$G$-crystal $(X,\gamma_X,\{e_i\}_{i\in I})$,
which is called the induced geometric $G$-crystals by 
unipotent $G$-crystal $({\bf X},f_X)$.
\end{thm2}
Due to the product structure of unipotent crystals, 
we can deduce a product structure of geometric crystals
derived from unipotent crystals,
which is a counterpart of tensor product structure
of Kashiwara's crystals.
We omit the explicit statement here
(see \cite{BK},\cite{N}).
\subsection{Crystals}

The notion ``crystal'' is introduced as a combinatorial
object by abstracting the properties of ``crystal bases'',
which has, in general, no corresponding $\uq$-module.

\begin{def2}
\label{crystal}
A {\it crystal} $B$ is a set endowed with the following maps:
\begin{eqnarray*}
&& wt:B\lar P,\\
&&\vep_i:B\lar\ZZ\sqcup\{-\infty\},\q
  \vp_i:B\lar\ZZ\sqcup\{-\infty\} \q{\hbox{for}}\q i\in I,\\
&&\eit:B\sqcup\{0\}\lar B\sqcup\{0\},
\q\fit:B\sqcup\{0\}\lar B\sqcup\{0\}\q{\hbox{for}}\q i\in I,\\
&&\eit(0)=\fit(0)=0,
\end{eqnarray*}
those maps satisfy the following axioms: for
 all $b,b_1,b_2 \in B$, we have
\begin{eqnarray}
&&\vp_i(b)=\vep_i(b)+\lan h_i,wt(b)\ran,\label{c1}\\
&&wt(\eit b)=wt(b)+\al_i{\hbox{ if  }}\eit b\in B,\label{c2}\\
&&wt(\fit b)=wt(b)-\al_i{\hbox{ if  }}\fit b\in B,\label{c3}\\
&&\eit b_2=b_1 \Longleftrightarrow \fit b_1=b_2\,\,(\,b_1,b_2 \in B),
\label{c4}\\
&&\vep_i(b)=-\ify
   \Longrightarrow \eit b=\fit b=0.\label{c5}
\end{eqnarray}
\end{def2}
The operators $\eit$ and $\fit$ are called the 
{\it Kashiwara operators}.
Indeed, if $(L,B)$ is a crystal base, then $B$ is a crystal.

\nd
{\sl Remark.} A {\it pre-crystal} is an object satisfying 
the conditions (\ref{c1})--(\ref{c3}).

Let us define $\til s_i:B\rightarrow B$ $(i\in I)$ 
(\cite{K3}) by
\[
 \til s_i(b)=
\begin{cases}
\eit^{-\lan wt(b),h_i\ran}(b)&\text{if }\lan wt(b),h_i\ran<0,\\
\fit^{\lan wt(b),h_i\ran}(b)&\text{if }\lan wt(b),h_i\ran\geq0.
\end{cases}
\]
Here note that we have $\til s_i^2={\rm id}_B$.
\begin{def2}
Let $B$ be a crystal.
\begin{enumerate}
\item
If the actions by $\{s_i\}_{i\in I}$ define
the action of the Weyl group $W$  on $B$, we call $B$ a $W$-{\it crystal}.
\item
If $\eit$ or $\fit$ is bijective, then we call $B$ a 
{\it free crystal}.
\end{enumerate}
\end{def2}
Note that if $B$ is a free crystal, then $\fit=\eit^{-1}$.
We frequently 
denote a free crystal $B$ by $(B,wt,\{\eit\}_{i\in I})$.

\subsection{Positive structure 
and Ultra-discretizations/Tropicalizations}

Let us recall the notions of 
positive structure and ultra-discretization/tropicalization.

The setting below is simpler than the ones in 
(\cite{BK},\cite{N} ), since it is sufficient for our 
purpose. 
Let $T=(\bbC^\times)^l$ be an algebraic torus over $\bbC$ and 
$X^*(T)\cong \ZZ^l$ (resp. $X_*(T)\cong \ZZ^l$) 
be the lattice of characters
(resp. co-characters)
of $T$. 
Set $R:=\bbC(c)$ and define
$$
\begin{array}{cccc}
v:&R\setminus\{0\}&\longrightarrow &\ZZ\\
&f(c)&\mapsto
&{\rm deg}(f(c)).
\end{array}
$$
Here note that for $f_1,f_2\in R\setminus\{0\}$, we have
\begin{equation}
v(f_1 f_2)=v(f_1)+v(f_2),\q
v\left(\frac{f_1}{f_2}\right)=v(f_1)-v(f_2)
\label{ff=f+f}
\end{equation}
Let $f=(f_1,\cd,f_n):T\rightarrow T'$ be a rational morphism 
between two algebraic tori $T=(\bbC^\times)^m$ and 
$T'=(\bbC^\times)^n$.
We define a map $\what f:X_*(T)\rightarrow X_*(T')$ by 
\[
(\what f(\xi))(c)
:=(c^{v(f_1(\xi(c))},\cd,c^{v(f_n(\xi(c)))}),
\]
where $\xi\in X_*(T)$.
Since $v$ satisfies (\ref{ff=f+f}), the map $\what f$
is an additive group homomorphism.
If we identify $X_*(T)$ (resp. $X_*(T')$)with $\ZZ^m$
(resp. $\ZZ^n$) by 
$\xi(c)=(c^{l_1},\cd,c^{l_m})\leftrightarrow 
(l_1,\cd,l_m)\in \ZZ^m$, 
we write 
\[
\what f(l_1,\cd,l_m)
:=(v(f_1(\xi(c))),\cd,v(f_n(\xi(c)))).
\]

A rational function $f(c)\in \bbC(c)$ $(f\ne0)$ is 
{\it positive} if $f$ can be expressed as a ratio
of polynomials with positive coefficients.

\nd
{\sl Remark.}
A rational function $f(c)\in \bbC(c)$ is positive 
if and only if $f(a)>0$ for any $a>0$
(pointed out by M.Kashiwara).

If $f_1,\,\,f_2\in R$ 
are positive, then we have (\ref{ff=f+f}) and 
\begin{equation}
v(f_1+f_2)={\rm max}(v(f_1),v(f_2)).
\label{+max}
\end{equation}
\begin{def2}[\cite{BK}]
Let 
$f=(f_1,\cd,f_n):T\rightarrow T'$ between
two algebraic tori $T,T'$ be a rational morphism as above. 
It is called {\it positive},
if the following two conditions are satisfied:
\begin{enumerate}
\item For any co-character 
$\xi:\bbC^\tm\rightarrow T$, the image of $\xi$
is contained in dom$(f)$.
\item For any co-character 
$\xi:\bbC^\tm\rightarrow T$, any $f_i(\xi(c))$ ($i\in I$)
is a positive rational function.
\end{enumerate}
\end{def2}

Denote by ${\rm Mor}^+(T,T')$ the set of 
positive rational morphisms from $T$ to $T'$.

\begin{lem2}[\cite{BK}]
\label{TTT}
For any positive rational morphisms 
$f\in {\rm Mor}^+(T_1,T_2)$             
and $g\in {\rm Mor}^+(T_2,T_3)$, 
 the composition $g\circ f$
is in ${\rm Mor}^+(T_1,T_3)$.
\end{lem2}

By Lemma \ref{TTT}, we can define a category ${\cal T}_+$
whose objects are algebraic tori over $\bbC$ and arrows
are positive rational morphisms.
\begin{lem2}[\cite{BK}]
For any algebraic tori $T_1$, $T_2$, $T_3$, 
and positive rational morphisms 
$f\in {\rm Mor}^+(T_1,T_2)$, 
$g\in {\rm Mor}^+(T_2,T_3)$, we have
$\what{g\circ f}=\what g\circ\what f.$
\end{lem2}
By this lemma, we obtain a functor 
\[
\begin{array}{cccc}
{\cal UD}:&{\cal T}_+&\longrightarrow &{{\hbox{\germ Set}}}\\
&T&\mapsto& X_*(T)\\
&(f:T\rightarrow T')&\mapsto& 
(\what f:X_*(T)\rightarrow X_*(T')))
\end{array}
\]


\begin{def2}[\cite{BK}]
Let $\chi=(X,\gamma,\{e_i\}_{i\in I})$ be a 
geometric crystal, $T'$ be an algebraic torus
and $\theta:T'\rightarrow X$ 
be a birational isomorphism.
The isomorphism $\theta$ is called 
{\it positive structure} on
$\chi$ if it satisfies
\begin{enumerate}
\item the rational morphism 
$\gamma\circ \theta:T'\rightarrow T$ is positive.
\item
For any $i\in I$, the rational morphism 
$e_{i,\theta}:\bbC^\tm \tm T'\rightarrow T'$ defined by
$e_{i,\theta}(c,t)
:=\theta^{-1}\circ e_i^c\circ \theta(t)$
is positive.
\end{enumerate}
\end{def2}
Let $\theta:T\rightarrow X$ be a positive structure on 
a geometric crystal $\chi=(X,\gamma,\{e_i\}_{i\in I}\})$.
Applying the functor ${\cal UD}$ 
to positive rational morphisms
$e_{i,\theta}:\bbC^\tm \tm T'\rightarrow T'$ and
$\gamma\circ \theta:T'\ra T$
(the notations are
as above), we obtain
\begin{eqnarray*}
\til e_i&:=&{\cal UD}(e_{i,\theta}):
\ZZ\tm X_*(T) \rightarrow X_*(T)\\
\til\gamma&:=&{\cal UD}(\gamma\circ\theta):
X_*(T')\rightarrow X_*(T).
\end{eqnarray*}
Now, for given positive structure $\theta:T'\rightarrow X$
on a geometric crystal 
$\chi=(X,\gamma,\{e_i\}_{i\in I})$, we associate 
the triplet $(X_*(T'),\til \gamma,\{\til e_i\}_{i\in I})$
with a free pre-crystal structure (see \cite[2.2]{BK}) 
and denote it by ${\cal UD}_{\theta,T'}(\chi)$.
By Lemma \ref{Verma}, we have the following theorem:

\begin{thm2}
For any geometric crystal 
$\chi=(X,\gamma,\{e_i\}_{i\in I})$ and positive structure
$\theta:T'\rightarrow X$, the associated pre-crystal 
${\cal UD}_{\theta,T'}(\chi)=
(X_*(T'),\til\gamma,\{\til e_i\}_{i\in I})$ 
is a free $W$-crystal (see \cite[2.2]{BK})
\end{thm2}

We call the functor $\cal UD$
{\it ``ultra-discretization''}
instead of ``tropicalization'' unlike in \cite{BK}.
And 
for a crystal $B$, if there
exists a geometric crystal $\chi$, 
an algebraic torus $T$ in ${\cal T}_+$ and a positive 
structure $\theta$ on $\chi$ such that 
${\cal UD}_{\theta,T}(\chi)\cong B$ as crystals, 
we call $\chi$
a {\it tropicalization} of $B$.

\renewcommand{\thesection}{\arabic{section}}
\section{Geometric crystals on unipotent groups}
\setcounter{equation}{0}
\renewcommand{\theequation}{\thesection.\arabic{equation}}
\newtheorem{pro4}{Proposition}[section]
\newtheorem{thm4}[pro4]{Theorem}
\newtheorem{lem4}[pro4]{Lemma}
\newtheorem{ex4}[pro4]{Example}
\newtheorem{cor4}[pro4]{Corollary}
\theoremstyle{definition}
\newtheorem{def4}[pro4]{Definition}

In this section, we associate a geometric/unipotent crystal structure
with unipotent subgroup $U^-$ of semi-simple algebraic group $G$.
In particular, for $G=SL_{n+1}(\bbC)$ we describe it explicitly.
\subsection{$U$-variety structure on $U^-$}
In this subsection, suppose that 
$G$ is a Kac-Moody group as in Sect.2.
As mentioned in \ref{uni-cry}, Borel subgroup $B^-$ has a 
$U$-variety structure. By the similar manner, we define
$U$-variety structure on $U^-$.
As in \ref{uni-cry}, 
the multiplication map $m$ in $G$ induces an open embedding;
$m:U^-\times B\hookrightarrow G,$
then this is a birational isomorphism. 
Let us denote the inverse birational isomorphism by $h$;
\[
 h:G\longrightarrow U^-\times B.
\]
Then we define the rational morphisms 
$\pi^{--}:G\rightarrow B^-$ and 
$\pi^+:G\rightarrow B$ by 
$\pi^{--}:={\rm proj}_{U^-}\circ h$ 
and $\pi^+:={\rm proj}_B\circ h$.
Now we define the rational $U$-action $\al_{U^-}$ on $U^-$ by 
\[
 \al_{U^-}:=\pi^{--}\circ m:U\times U^-\longrightarrow U^-,
\]
Then we obtain 
\begin{lem4}
\label{XU}
A pair ${\bf U}^-=(U^-,\al_{U^-})$ is a 
$U$-variety on a unipotent subgroup $U^-\subset G$.
\end{lem4}
\subsection{Unipotent/Geometric crystal structure on $U^-$}
In order to define a unipotent crystal structure on $U^-$,
let us construct a $U$-morphism 
$F:U^-\rightarrow B^-$.

The multiplication map $m$ in $G$ induces an open embedding;
$m:U^-\times T\times U\hookrightarrow G,$
which is a birational isomorphism. 
Thus, by the similar way as above, we obtain 
the rational morphism
$\pi^{0}:G\rightarrow T$.
Here note that we have 
\begin{equation}
\pi^-(x)=\pi^{--}(x)\pi^0(x)\q(x\in G).
\label{--0}
\end{equation}
Now, we give a sufficient condition for existence of 
$U$-morphism $F$.
\begin{lem4}
\label{suf}
Let ${\cal T}:U^-\rightarrow T$ be a rational morphism satisfying:
\begin{equation}
{\cal T}(\pi^{--}(xu))=\pi^0(xu){\cal T}(u),
\q\text{for $x\in U$ and $u\in U^-$. }
\label{T}
\end{equation}
Defining a morphism $F:U^-\rightarrow B^-$ by
\begin{equation}
\begin{array}{cccc}
F:&U^-&\longrightarrow &B^-\\
&u&\mapsto& u{\cal T}(u),
\end{array}
\end{equation}
then the morphism $F$ is a $U$-morphism $U^-\rightarrow B^-$.
\end{lem4}
{\sl Proof.}
We may show 
\begin{equation}
F(\al_{U^-}(x,u))=\al_{B^-}(x,F(u)),
\q\text{for $x\in U$ and $u\in U^-$}.
\label{FU}
\end{equation}
As for the left-hand side of (\ref{FU}), we have
\[
F(\al_{U^-}(x,u))= \pi^{--}(xu){\cal T}(\pi^{--}(xu))
=\pi^{--}(xu)\pi^0(xu){\cal T}(u),
\]
where the last equality is due to (\ref{T}).
On the other hand, the right-hand side
 of (\ref{FU}) is written by:
\[
 \al_{B^-}(x,F(u))=\pi^-(xu{\cal T}(u))
=\pi^{--}(xu{\cal T}(u))\pi^0(xu{\cal T}(u))
=\pi^{--}(xu)\pi^0(xu){\cal T}(u)
\]
where the second equality is due to (\ref{--0}) and the 
third equality is obtained by the fact that 
${\cal T}(u)\in T\subset B$.
Now we get (\ref{FU}).\qed

Let us verify that there exists such $U$-morphism $F$
or rational morphism ${\cal T}$ for semisimple cases.
Suppose that $G$ is semisimple in the rest of this section.

Let $\Lm_i\in P_+$ ($i=1,\cd,n$)
be a fundamental weight and $L(\Lm_i)$
be a corresponding irreducible highest weight 
$\ge$-module, where $\ge$ is a complex semi-simple
Lie algebra associated with $G$.
Let $L^*(\Lm_i)$ be a contragredient module of 
$L(\Lm_i)$ as in Sect.2
and  fix a  highest (resp. lowest) weight vector $u_{\Lm_i}$
(resp. $v_{\Lm_i}$) and $v^*_{\Lm_i}\in L^*(\Lm_i)$ be the same vector as 
$v_{\Lm_i}$
such that $\lan v_{\Lm_i},v^*_{\Lm_i}\ran=1$.
Now, let us define 
a function $f_i:U^-\rightarrow \bbC$ $(i\in I)$ as a matrix
coefficient: 
\begin{equation}
f_i(g)=\lan g\cdot u_{\Lm_i},v^*_{\Lm_i}\ran.
\label{fi}
\end{equation}
We define a rational morphism ${\cal T}:U^-\rightarrow T$ by
\begin{equation}
{\cal T}(u):=\prod_{i\in I}\al_i^\vee(f_i(u)^{-1}).
\label{TT}
\end{equation}
and define a morphism $F:U^-\rightarrow B^-$ by
\begin{equation}
F(u):=u\cdot \prod_{i\in I}\al_i^\vee(f_i(u)^{-1}).
\label{FF}
\end{equation}
\begin{lem4}
The morphism $F:U^-\rightarrow B^-$  is a $U$-morphism.
\end{lem4}
{\sl Proof.}
Let us verify that $\cal T$ satisfies (\ref{T}).
For $x\in U$ and $u\in U^-$ such that 
$xu\in Im(U^-\times T\times U\hookrightarrow G)$,
let $u^-\in U^-$, $u^0\in T$ and $u^+\in U$  
be the unique elements satisfying $u^-u^0u^+=xu$, 
{\it i.e.,} 
$\pi^{--}(xu)=u^-$, $\pi^0(xu)=u^0$ and $\pi(xu)=u^+$.
Since $\lan\q,\q\ran$ is a contragredient bilinear form
and the fact that $g\cdot v^*_{\Lm_i}=v^*_{\Lm_i}$ for 
any $g\in U^-$, we have 
\begin{equation}
\lan xu\cdot u_{\Lm_i}, v^*_{\Lm_i}\ran
=\lan u\cdot u_{\Lm_i}, \omega(x)\cdot v^*_{\Lm_i}\ran
=\lan u\cdot u_{\Lm_i}, v^*_{\Lm_i}\ran.
\label{xu*}
\end{equation}
On the other hand, since $g\cdot u_{\Lm_i}=u_{\Lm_i}$ for $g\in U$,
we have 
\begin{eqnarray}
&& \lan xu\cdot u_{\Lm_i}, v^*_{\Lm_i}\ran
=\lan \pi^-(xu)\pi^0(xu)\pi(xu)\cdot u_{\Lm_i}, v^*_{\Lm_i}\ran
\nn\\
&&=\lan \pi^-(xu)\pi^0(xu)\cdot u_{\Lm_i}, v^*_{\Lm_i}\ran
=\Lm_i(\pi^0(xu))\lan \pi^-(xu)\cdot u_{\Lm_i}, v^*_{\Lm_i}\ran,
\label{xu0}
\end{eqnarray}
where $\Lm_i\in X^*(T)$ such that 
$\Lm_i(\al_j^\vee(c))=c^{\del_{i,j}}$.
Hence, by (\ref{xu*}), (\ref{xu0}), we have
\begin{eqnarray*}
f_i(\pi^-(xu))&=&\lan \pi^-(xu)\cdot u_{\Lm_i}, v^*_{\Lm_i}\ran
=\Lm_i(\pi^0(xu))\lan \pi^-(xu)\cdot u_{\Lm_i}, v^*_{\Lm_i}\ran\\
&=&\Lm_i(\pi^0(xu))^{-1}\lan uu_{\Lm_i}, v^*_{\Lm_i}\ran
=\Lm_i(\pi^0(xu))^{-1}f_i(u).
\end{eqnarray*}
By the formula
\[
 \prod_i\al_i^\vee(\Lm_i(t))=t, \q(t\in T),
\]
and the definitions of $\cal T$ and $F$, 
we obtained (\ref{T}).\qed

\begin{cor4}
Suppose that $G$ is semi-simple. 
Then $(U^-,F)$ is a unipotent crystal.
\end{cor4}
As we have seen in \ref{u->g}, we can associate 
geometric crystal structure with the unipotent subgroup $U^-$
since it has a unipotent crystal structure.

It is trivial that the function 
$\vp_i:U^-\rightarrow \bbC$ is not identically zero.
Thus, defining the morphisms 
$e_i:\bbC^\times\times U^-\rightarrow U^-$ and 
$\gamma_{U^-}:U^-\rightarrow T$ by
\begin{equation}
e_i(c,u)=e_i^c(u):=x_i(\frac{c-1}{\vp_i(u)})(u),\qq
\gamma_{U^-}(u):={\cal T}(u),\qq
\text{($u\in U^-$ and $c\in \bbC^\times$)},
\end{equation}
It follows from Theorem \ref{U-G}:
\begin{thm4}
If $G$ is semi-simple, then the triplet 
$\chi_{U^-}:=(U^-,\gamma_{U^-},\{e_i\}_{i\in I})$ is a geometric 
crystal.
\end{thm4}
\section{$SL_{n+1}(\bbC)$-case}
\setcounter{equation}{0}

We see the result of the previous section in the 
$SL_{n+1}(\bbC)$-case more explicitly. 

We identify unipotent subgroup $U^-$ with the set
of lower triangular matrices 
whose diagonal part is an identity matrix.

First, let us describe
 the morphism $\cal T:U^-\rightarrow T$.
For $i\in I:=\{1,\cd,n\}$ 
and $u=(a_{ij})_{1\leq i,j\leq n+1}\in U^-$, 
let $u^{(i)}$ be the submatrix with size $i$ as:
\[
 u^{(i)}:=(a_{i,j})_{n-i+2\leq i\leq n+1, 1\leq j\leq i},
\]
{\it i.e.,}
\[
u=\left(
\begin{array}{ccccc}
1&&&&\\
&1&&{\hbox{\Huge0}}&\\
&&\cdot&&\\
&&&\cdot&\\
&&&&1
\end{array}
\right)\in U^-
\]
\vskip-1.5cm
\hspace{4.9cm}$\underbrace{\Fsquare(10mm,u^{(i)})}_i$

\nd
and set 
$ m_i(u):=\det(u^{(i)}).$

Let $V=\bbC^{n+1}$ be the $n+1$-dimensional vector space
with the basis $\{u_1,u_2,\cd,u_{n+1}\}$. 
We can identify $V$ with 
the vector representation $L(\Lm_1)$ of $\ssl_{n+1}$
by the standard way.
Indeed, the explicit actions are given by:
\[
e_i(u_j)=\del_{i+1,j}u_{i-1}\q f_i(u_j)=\del_{i,j}u_{i+1},
\]
where $e_i=E_{i,i+1}$ and $f_i=E_{i+1,i}$(matrix unit), and
we set $e_i(u_1)=0$ and $f_i(u_{n+1})=0$  for all $i\in I$,
which implies that $u_1$ is the highest weight vector and
$u_{n+1}$ is the lowest weight vector.
Then 
we have the isomorphism between the fundamental representation
$L(\Lm_k)$ $(1\leq k\leq n)$ 
and the $k$-th anti-symmetric tensor module
$\bigwedge^k(V)$. 
Let us fix  
\begin{equation}
 u_{\Lm_k}:=u_1\wedge u_2\wedge\cd\wedge u_{k}
\text{(resp. $
 v_{\Lm_k}:=u_{n-k+2}\wedge u_{n-k+3}
\wedge\cd\wedge u_{n+1}$)}
\label{h-l}
\end{equation}
the highest (resp. lowest)weight vector in 
$L(\Lm_k)\cong \bigwedge^k(V)$.
In this setting, we have
\begin{lem4}
$f_i\equiv m_i$ on $U^-$ for all $i=1,\cd,n$.
\end{lem4}
{\sl Proof.}
For $g=(g_{ij})\in U^-$, we have
\[
g\cdot u_i=u_i+\sum_{i<j}g_{ji}u_j.
\]
Let us see the coefficient of the vector 
$v_{\Lm_k}:=u_{n-k+2}\wedge u_{n-k+3}\wedge\cd\wedge u_{n+1}$
in $g\cdot u_{\Lm_k}$. We have 
\begin{eqnarray*}
g\cdot u_{\Lm_k}&=&g\cdot u_1\wedge g\cdot u_2\wedge
\cd\wedge g\cdot u_{k}\\
&=& 
(u_1+\sum_{1<j}g_{j1}u_j)\wedge \cd\wedge
(u_k+\sum_{k<j}g_{jk}u_j).
\end{eqnarray*}
Thus, the coefficient of the vector 
$u_{j_1}\wedge\cd\wedge u_{j_k}$ is 
$g_{j_1\,1}\cd g_{j_k\,k}$. Hence, we obtain the coefficient 
of $v_{\Lm_k}$ as
\[
\sum_{\sigma\in \mathfrak{S}_{k}}
{\rm sgn}(\sigma)g_{n-\sigma(1)+2\,1} \cd g_{n-\sigma(k)+2\,k}
=m_i(u),
\]
by using $v_{\sigma(1)}\wedge 
v_{\sigma(2)}\wedge\cd\wedge v_{\sigma(k)}
={\rm sgn}(\sigma)v_1\wedge v_2\wedge \cd\wedge v_k$.
On the other hand,  the coefficient of the lowest weight vector
gives the function $f_i(u)$. Then we get the desired result.\qed

Next, let us see the action of $e_i^c$ on $U^-$.
Indeed, the action of $e_i^c$ is described simply by;
\[
 e_i^\al(u)=x_i(\frac{\vp_i(u)}{\al-1})\cdot u
\cdot x_i(\frac{1-\al}{\al\vp_i(u)})\cdot\al_i^{-1}(\al).
\]
Here, for the later purpose, we consider the following
subset $B^u$ of $U^-$ and describe the action of $e_i^\al$ on it:
\begin{equation}
B^u:=
\left\{
Y(a)=\begin{array}{l}
y_n(a_{1,n})y_{n-1}(a_{1,n-1})\cd y_1(a_{1,1})\times \\
\times y_n(a_{2,n})\cd\cd y_2(a_{2,2})\times \\
\cd\cd\cd\\
\times y_n(a_{n,n}) \\
\end{array}
:a_{i,j}\in\mathbb{C}^\times
\right\}
\subset U^-.
\end{equation}
It is easy to see that $B^u$ is an open dense subset in $U^-$
and isomorphic to the algebraic torus 
$T_0:=(\bbC^\times)^{\frac{n(n+1)}{2}}$ by:
\begin{eqnarray*}
T_0=(\bbC^\times)^{\frac{n(n+1)}{2}}&\mapright{\sim} & B^u,\\
a=(a_{i,j})_{1\leq i\leq j\leq n}&\mapsto&Y(a),
\end{eqnarray*}
which gives a birational isomorphism 
$\theta:T_0=(\bbC^\times)^{\frac{n(n+1)}{2}}
\rightarrow B^u\hookrightarrow U^-$.

Furthermore, we have
\begin{lem4}
\label{phi-fi}
$\vp_i(Y(a))=\sum_{k=1}^{i}a_{k,i},\q
f_i(Y(a))=\prod_{k=1}^{i}\prod_{j=k}^{n-i+k}a_{k,j}.$
\end{lem4}
{\sl Proof.}
Since for $u\in U^-$, $\vp_i(u)$ is given as a 
$(i+1,i)$-entry and the $(i+1,i)$-entry of $Y(a)$ is 
$\sum_{k=1}^{i}a_{k,i}$, we obtained the first result.

For a word $\io=i_1,\cd,i_m$, we set
$f_\io:=f_{i_1}\cd f_{i_m}$. 
For a fixed reduced longest word 
$\io_0=n,\cd,1,n\cd,2,n\cd,3,\cd,n,n-1,n$,
there exists the unique subword 
\[
\io_i:=\underbrace{n-i+1,\cd,1},
\underbrace{n-i+2,\cd,2},\cd,
\underbrace{n-1,\cd,i-1},\underbrace{n,\cd,i},
\]
such that 
\begin{equation}
f_{\io_i}u_{\Lm_i}=v_{\Lm_i},
\label{f-lm}
\end{equation}
where $u_{\Lm_i}$ (resp. $v_{\Lm_i}$) is the highest 
(resp. lowest) weight vector
of $L(\Lm_i)$ as in (\ref{h-l}).
Since $f_i^2(L(\Lm_k))=\{0\}$, 
we have $y_i(a)=1+af_i$ on $L(\Lm_k)$  and then
\begin{eqnarray*}
Y(a)&=&(1+a_{1,n}f_n)\cd(1+a_{1,1}f_1)\cd(1+a_{n,n}f_n)\\
&=& \sum_{\io:\text{subword of }\io_0}a_{\io}f_\io,
\end{eqnarray*}
where $a_\io$ is a coefficient of $f_\io$ and $a_{\io}f_\io=1$
if $\io$ is empty. Hence by (\ref{f-lm}), we have
\[
f_i(Y(a))= \lan Y(a)u_{\Lm_i},v^*_{\Lm_i}\ran=a_{\io_i}=
\prod_{k=1}^{i}\prod_{j=k}^{n-i+k}a_{k,j}.
\]

\qed

Let us see the rational action 
$e_i^\al:B^u\rightarrow B^u$. 
\begin{pro4}
We have 
$e_i^\al Y((a_{k,j})_{1\leq k\leq j\leq n})
=Y((a'_{k,j})_{1\leq k\leq j\leq n}),
$ 
\begin{equation}
a'_{k,j}:=
\begin{cases}
C^{(i)}_{k}a_{k,i-1}&\text{if $j=i-1$},\\
\frac{a_{k,i}}{C^{(i)}_{k-1}C^{(i)}_{k}}&\text{if $j=i$},\\
C^{(i)}_{k-1}a_{k,i+1}&\text{if $j=i+1$},\\
a_{k,j}&\text{otherwise},
\end{cases}
\end{equation}
where
\[
 C^{(i)}_k:=\frac{\al(a_{1,i}+\cd+a_{k,i})+a_{k+1,i}+\cd+a_{i,i}}
{a_{1,i}+\cd+a_{i,i}}\q(1\leq k\leq i\leq n).
\]
\end{pro4}
{\sl Proof.}
We recall the formula:
\begin{eqnarray}
&&x_i(a)y_j(b)=
\left\{
\begin{array}{ll}
 y_i(\frac{b}{1+ab})\al_i^\vee(1+ab)x_i(\frac{a}{1+ab}) &
{\rm if }\,\,i=j\\
 y_j(b)x_i(a)&{\rm if}\,\,i\ne j
\end{array}
\right.,
\label{xy=yax}\\
&&\al_i^\vee(a)x_j(b)=x_j(a^{a_{ij}}b)
\al_i^\vee(a),
\q
\al_i^\vee(a)y_j(b)=y_j(a^{-a_{ij}}b)
\al_i^\vee(a).
\label{ae}
\end{eqnarray}
Using these formula repeatedly, we have
\[
 x_i(c)\cdot Y(a)=Y(a')\cdot \al_i^\vee(1+c\vp_i(Y(a)))\cdot
x_i(\frac{c}{1+c\vp_i(Y(a))}),
\]
where $c=(\al-1)/(a_{1,i}+\cd+a_{i,i})$ and 
$\vp_i(Y(a))=a_{1,i}+\cd+a_{i,i}$. \qed

Now, we consider the following birational isomorphism:
\begin{eqnarray*}
 \xi:T_0=(\bbC^\times)^{\frac{n(n+1)}{2}}&\longrightarrow &
T_0=(\bbC^\times)^{\frac{n(n+1)}{2}},\\
(a_{i,j})_{1\leq i\leq j\leq n}&\mapsto&
(A_{i,j})_{1\leq i\leq j\leq n}
\end{eqnarray*}
where
\[
A_{i,j}:=
\frac{a_{i,j}a_{i-1,j-1}\cd a_{1,j-i+1}}
{a_{i-1,j}a_{i-2,j-1}\cd a_{1,j-i+2}}\quad
(1\leq i\leq j\leq n).
\]
The inverse morphism is  given by
\begin{equation}
 a_{i,j}:=
\frac{A_{i,j}A_{i-1,j}\cd A_{1,j}}
{A_{i-1,j-1}A_{i-2,j-1}\cd A_{1,j-1}}\quad
(1\leq i\leq j\leq n).
\label{a-A}
\end{equation}
We can describe explicitly 
\begin{eqnarray*}
\xi\circ e_i^\al\circ \xi^{-1}:(\bbC^\times)^{\frac{n(n+1)}{2}}
&\longrightarrow &(\bbC^\times)^{\frac{n(n+1)}{2}},\\
(A_{k,j})_{1\leq k\leq j\leq n}
&\mapsto&
(A'_{k,j})_{1\leq k\leq j\leq n},
\end{eqnarray*}
where 
\[
 A'_{k,j}=
\begin{cases}
A_{k,j}&\text{if $j\ne i,i-1$}\\
\al^{(i)}_k\cdot A_{k,i-1}&\text{if $j=i-1$}\\
({\al^{(i)}_k})^{-1}\cdot A_{k,i}&\text{if $j=i$}
\end{cases}
\]
\begin{equation}
\al^{(i)}_k=
\frac{\displaystyle
\al\sum_{1\leq j\leq k}
{\frac{\prod_{l=1}^j A_{l,i}}
{\prod_{l=1}^{j-1} A_{l,i-1}}}
+\sum_{k<j\leq i}
{\frac{\prod_{l=1}^j A_{l,i}}
{\prod_{l=1}^{j-1} A_{l,i-1}}}}
{\displaystyle
\al\sum_{
1\leq j\leq k-1}
{\frac{\prod_{l=1}^j A_{l,i}}
{\prod_{l=1}^{j-1} A_{l,i-1}}}
+\sum_{j=k}^{i}
{\frac{\prod_{l=1}^j A_{l,i}}
{\prod_{l=1}^{j-1} A_{l,i-1}}}}.
\label{alik}
\end{equation}
Set $\hat\theta:=\theta\circ\xi^{-1}:
(\bbC^\times)^{\frac{n(n+1)}{2}}
\mapright{\xi^{-1}}
(\bbC^\times)^{\frac{n(n+1)}{2}}
\mapright{\theta} U^-$.
\begin{thm4}
The morphism $\hat\theta$ gives a  positive structure 
on the geometric crystal $\chi_{U^-}$.
\end{thm4}
{\sl Proof.}
The explicit form of 
\begin{eqnarray*}
\hat\theta^{-1}\circ e_i^{\al}\circ\hat\theta:
(\bbC^\times)^{\frac{n(n+1)}{2}}\times \bbC^\times
&\longrightarrow& (\bbC^\times)^{\frac{n(n+1)}{2}}\\
(A_{k,j})_{1\leq k\leq j\leq n}&\mapsto &
(A_{k,j}')_{1\leq k\leq j\leq n}
\end{eqnarray*}
is given as (\ref{alik}), which is trivially positive.
Then let us show the positivity of $\gamma_{U^-}\circ\hat\theta$.
For $Y(a)\in B^u$, we have 
$\gamma_{U^-}(Y(a))=\prod_i\al_i^\vee(f_i(Y(a))^{-1})$ and 
by Lemma \ref{phi-fi} 
the explicit form of $f_i(Y(a))$ is given.
Substituting (\ref{a-A}) in it, we obtain
\begin{equation}
\gamma_{U^-}\circ\hat\theta((A_{k,j})_{1\leq k\leq j\leq n})
=\gamma_{U^-}\circ\theta\circ\xi^{-1}
((A_{k,j})_{1\leq k\leq j\leq n})
=\prod_{i=1}^{n}\al^\vee_i(\prod_{\begin{array}{c}
\scriptstyle 1\leq k\leq i\\ \scriptstyle
i\leq j\leq n\end{array}}A_{k,j})^{-1},
\label{gamma-theta}
\end{equation}
which implies that $\gamma\circ\hat\theta$ is positive.\qed
\section{Tropicalization of 
Geometric Crystals on $U^-$ and generalized Young Tableaux}
\setcounter{equation}{0}
\subsection{Crystal structure on Young tableaux}
Let us recall the crystal structure on Young tableaux
where the terminology ``Young tableaux'' means ``semi-standard 
tableaux'' in \cite{KN}.
For a partition $\lm=(\lm_1,\lm_2,\cd,\lm_n)$, set 
\[
B(\lm)=\{\text{ Young tableau  of shape }\lm
\text{ with contents }1,2,\cd,n+1\},
\]
which gives an $A_n$-crystal of irreducible highest weight 
$U_q(\ssl_{n+1})$-module $V(\lm)$\cite{KN}.

In order to describe the action of $\til e_i^\beta(b)$ 
$(\beta\geq0)$ explicitly, let us recall how to
construct $B(\lm)$ following \cite{KN},\cite{K5}.

\def\bx(#1){\fsquare(5mm,#1)}
\def\bbx{B_{\fsquare(2mm,)}}
Let $V_{\fsquare(2mm,)}:=
V(\Lm_1)$ be the vector representation
of $U_q(\ssl_{n+1})$, which is the irreducible 
highest weight module with the highest weight $\Lm_1$ and 
let 
\[
 \bbx:=\left\{\bx(1),\bx(2),\cd,\bx(\tiny n+1)\right\}
\]
be the crystal of $V_{\fsquare(2mm,)}$.
The explicit actions of $\til e_i$ and $\til f_i$ are given as
follows(\cite{KN}):
\[
 \til e_i\bx(j)=\del_{i+1,j}\bx(i-1),\qq
 \til f_i\bx(j)=\del_{i,j}\bx(i+1).
\]
We realize $B(\lm)$ by embedding into
$\bbx^N$ ($N=|\lm|$), which follows the way of 
embedding $V(\lm)\hookrightarrow V_{\fsquare(2mm,)}^N$. 
In \cite{KN},
the ``{\it Japanese reading}'' is introduced, which gives the 
embedding by reading entries in a Young tableau
column by column. 
But here we take so-called ``{\it arabic reading}''
\cite{K5}, which gives the 
embedding by reading entries in a Young 
tableau row by row from right to left since it matches
what we do below.
\begin{ex4}
\begin{enumerate}
\item {\it Japanese reading}
\[
 \fourtwoonebox(a,b,c,d,e,f,g)=\left(\bx(d)\right)
\ot\left(\bx(c)\right)\ot
\left(\bx(b)\ot\bx(f)\right)
\ot\left(\bx(a)\ot\bx(e)\ot\bx(g)\right)\in \bbx^{\ot 7}
\]
\item
{\it Arabic reading}
\[
 \fourtwoonebox(a,b,c,d,e,f,g)
=\left(\bx(d)\ot\bx(c)\ot\bx(b)\ot\bx(a)\right)
\ot\left(\bx(f)\ot\bx(e)\right)\ot\left(\bx(g)\right)
\in \bbx^{\ot 7}
\]
\end{enumerate}
\end{ex4}
The description of the actions of 
$\eit$ and $\fit$ on $B(\lm)$ in \cite{KN} is as follows:
Let $\{(+),(-)\}$ (resp. $\{(0)\}$)be the crystal of 
the irreducible $U_q(\ssl_2)$-module $V_{\fsquare(2mm,)}$
(resp. $V(0)$). If we consider the actions of $\eit$ and $\fit$
on a tensor product $\bbx^N$, we can identify (\cite{KN},2.1.),
\begin{equation}
\bx(i)=(+), \q \bx(i+1)=(-),\q \bx(j)=(0)\,\,(j\ne i,i+1).
\label{+-0}
\end{equation}
Let $b\in B(\lm)$ be in the following form:

\begin{equation}
\hspace{-1.5cm}
\unitlength 0.1in
\begin{picture}( 64.2500, 29.9000)(  9.3000,-34.8000)
%
\special{pn 13}%
\special{pa 1326 490}%
\special{pa 1326 3480}%
\special{fp}%
%
\special{pn 13}%
\special{pa 1306 2290}%
\special{pa 4136 2300}%
\special{fp}%
%
\special{pn 13}%
\special{pa 1326 2680}%
\special{pa 4136 2680}%
\special{fp}%
%
\special{pn 13}%
\special{pa 1316 3090}%
\special{pa 1316 3090}%
\special{fp}%
%
\special{pn 13}%
\special{pa 1326 3090}%
\special{pa 2926 3080}%
\special{fp}%
%
\special{pn 13}%
\special{pa 2776 2670}%
\special{pa 2776 3080}%
\special{fp}%
\put(20.9500,-29.0000){\makebox(0,0){$\mathbf B_{i+1,i+1}$}}%
%
\special{pn 13}%
\special{pa 2886 2300}%
\special{pa 2886 2670}%
\special{fp}%
\put(20.7500,-25.0000){\makebox(0,0){$\mathbf B_{i,i}$}}%
\put(36.1500,-24.9000){\makebox(0,0){$\mathbf B_{i,i+1}$}}%
\put(11.3500,-24.9000){\makebox(0,0){$i$}}%
\put(11.5500,-29.1000){\makebox(0,0){$i+1$}}%
%
\special{pn 13}%
\special{pa 1316 490}%
\special{pa 7356 490}%
\special{fp}%
%
\special{pn 13}%
\special{pa 2926 880}%
\special{pa 7316 880}%
\special{fp}%
%
\special{pn 13}%
\special{pa 2516 1280}%
\special{pa 6526 1280}%
\special{fp}%
%
\special{pn 13}%
\special{pa 6926 490}%
\special{pa 6926 870}%
\special{fp}%
%
\special{pn 13}%
\special{pa 5636 490}%
\special{pa 5636 870}%
\special{fp}%
%
\special{pn 13}%
\special{pa 4316 490}%
\special{pa 4316 870}%
\special{fp}%
%
\special{pn 13}%
\special{pa 5396 880}%
\special{pa 5396 1260}%
\special{fp}%
%
\special{pn 13}%
\special{pa 4196 890}%
\special{pa 4196 1270}%
\special{fp}%
%
\special{pn 13}%
\special{pa 2916 890}%
\special{pa 2916 1270}%
\special{fp}%
\put(34.5500,-11.0000){\makebox(0,0){$\mathbf B_{2,i}$}}%
\put(49.8500,-10.8000){\makebox(0,0){$\mathbf B_{2,i+1}$}}%
\put(50.0500,-6.8000){\makebox(0,0){$\mathbf B_{1,i}$}}%
\put(62.6500,-6.8000){\makebox(0,0){$\mathbf B_{1,i+1}$}}%
%
\special{pn 13}%
\special{pa 1326 880}%
\special{pa 2906 880}%
\special{fp}%
%
\special{pn 13}%
\special{pa 1316 1280}%
\special{pa 2546 1280}%
\special{fp}%
%
\special{pn 13}%
\special{pa 1316 1670}%
\special{pa 5696 1670}%
\special{fp}%
\put(11.0500,-7.6000){\makebox(0,0)[lb]{1}}%
\put(11.1500,-11.9000){\makebox(0,0)[lb]{2}}%
\end{picture}%
\label{YT}
\end{equation}
where $B_{i,j}:=
\sharp\{j\hbox{ in the $i$-th row }\}$.
If we consider the actions of $\eit$ and $\fit$,
by the ``arabic reading'' and 
(\ref{+-0}) we can identify:
\[
 b=v_1\ot \cd\ot v_{i+1},
\]
where 
\begin{equation}
 v_k:=(-)^{\ot B_{k,i+1}}\ot (+)^{\ot B_{k,i}}\q
(1\leq k\leq i),\q
v_{i+1}=(-)^{\ot B_{i+1,i+1}}
\label{vk}
\end{equation}
For any $i\in I$ and $\beta\in \ZZ_{\geq0}$ there exist unique
$\beta^{(i)}_k\in \ZZ_{\geq 0}$ 
$(1\leq k\leq i+1)$ such that 
\begin{equation}
\til e_i^{\beta}(v_1\ot \cd\ot v_{i+1})
=\til e_i^{\beta^{(i)}_1}(v_1)\ot\cd\ot
\til e_i^{\beta^{(i)}_{i+1}}(v_{i+1}),
\label{ebeta}
\end{equation}
and $\beta=\sum_{1\leq k\leq i+1}\beta^{(i)}_k$.
Note that on each component, we have
\[
\til e_i^{\beta^{(i)}_k}(v_k):=
(-)^{\ot( B_{k,i+1}-\beta^{(i)}_k)}
\ot (+)^{\ot (B_{k,i}+\beta^{(i)}_k)}.
\]
Let us see the explicit form of $\beta^{(i)}_k$, in 
order to describe the action of $\eit^\beta$ on $b$.
For the purpose, we prepare the following formula:
\begin{lem4}[\cite{K5}]
\label{lem-ten}
Let $B_1$, $B_2,\cd,B_l$ be crystals. 
For $v_k\in B_k$ and $i\in I$, 
set $b_k:=\vep_i(v_k)-\sum_{1\leq j<k}\lan h_i,wt(v_j)\ran$.
Then, we have
\[
 \eit^c(v_1\ot\cd\ot v_l)=
\til e_i^{c_1}(v_1)\ot\cd\ot
\til e_i^{c_{l}}(v_{l}), 
\]
where
\begin{equation}
c_k=
{\rm max}(c+\mathop{\hbox{\rm max}}_{1\leq j\leq k}(b_j),
\mathop{\hbox{\rm max}}_{k<j\leq l}(b_j))
-{\rm max}(c+\mathop{\hbox{\rm max}}_{1\leq j<k}(b_j),
\mathop{\hbox{\rm max}}_{k\leq j\leq l}(b_j)).
\end{equation}
\end{lem4}
Applying this lemma to (\ref{ebeta}), we have
\begin{pro4}
\label{ex-ebeta}
Under the setting (\ref{vk}) and (\ref{ebeta}), 
\begin{eqnarray*}
\beta^{(i)}_k&=&{\hbox{\rm max}}\left(
\beta+\mathop{\hbox{\rm max}}_{
\scriptstyle
1\leq j\leq k}
\left(\sum_{l=1}^{j}
B_{l,i+1}-
\sum_{l=1}^{j-1}
B_{l,i}\right),
\mathop{\hbox{\rm max}}_{
\scriptstyle
k<j\leq i}
\left({\sum_{l=1}^{j}
B_{l,i+1}-
\sum_{l=1}^{j-1}
B_{l,i}}\right)
\right) 
\label{betaik}  \\
&&
-{\hbox{\rm max}}\left(
\beta+\mathop{\hbox{\rm max}}_{
\scriptstyle
1\leq j<k}
\left({\sum_{l=1}^{j}
B_{l,i+1}-
\sum_{l=1}^{j-1}
B_{l,i}}\right),
\mathop{\hbox{\rm max}}_{
\scriptstyle
k\leq j\leq i}
\left({\sum_{l=1}^{j}
B_{l,i+1}-
\sum_{l=1}^{j-1}
B_{l,i}}\right)
\right)\q(1\leq k\leq i),\\
\beta^{(i)}_{i+1}&=&0
\end{eqnarray*}
\end{pro4}
{\sl Proof.}
By (\ref{vk}), we have
$\vep_i(v_k)=B_{k,i+1},$ 
$\lan h_i,wt(v_j)\ran=B_{j,i}-B_{j,i+1}$.
Applying this to Lemma \ref{lem-ten}, we obtain
\begin{eqnarray}
&&
\beta^{(i)}_k={\hbox{\rm max}}\left(
\beta+\mathop{\hbox{\rm max}}_{
\scriptstyle
1\leq j\leq k}
\left(\sum_{l=1}^{j}
B_{l,i+1}-
\sum_{l=1}^{j-1}
B_{l,i}\right),
\mathop{\hbox{\rm max}}_{
\scriptstyle
k<j\leq i+1}
\left({\sum_{l=1}^{j}
B_{l,i+1}-
\sum_{l=1}^{j-1}
B_{l,i}}\right)
\right)  \nn \\
&&
-{\hbox{\rm max}}\left(
\beta+\mathop{\hbox{\rm max}}_{
\scriptstyle
1\leq j<k}
\left({\sum_{l=1}^{j}
B_{l,i+1}-
\sum_{l=1}^{j-1}
B_{l,i}}\right),
\mathop{\hbox{\rm max}}_{
\scriptstyle
k\leq j\leq i+1}
\left({\sum_{l=1}^{j}
B_{l,i+1}-
\sum_{l=1}^{j-1}
B_{l,i}}\right)
\right).\label{bk2}
\end{eqnarray}
Since $B_{i+1,i+1}\leq B_{i,i}$, we have
\[
\sum_{l=1}^{i}
B_{l,i+1}-
\sum_{l=1}^{i-1}
B_{l,i}\geq
\sum_{l=1}^{i+1}
B_{l,i+1}-
\sum_{l=1}^{i}
B_{l,i}.
\]
Hence, we can neglect $j=i+1$ in the formula (\ref{bk2}).
\qed

\nd
{\sl Remark.}
The formula $\beta^{(i)}_k$ does not depend on $B_{i,i}$
or $B_{i+1,i+1}$.
\subsection{Generalized Young tableaux
and its crystal structure}

Let $b\in B(\lm)$ be a Young tableau as in (\ref{YT}).
The $B_{i,j}$'s have several constraints, {\it e.g.,}
\[
 B_{i,j}\geq0, \q
\sum_{i\leq j\leq k}B_{i,j}\geq 
\sum_{i+1\leq j\leq k+1}B_{i+1,j},
\]
which come from the conditions for being Young tableaux.

Now, forgetting such constraints on $B_{i,j}$'s,
we obtain a free $\ZZ$-lattice $B^\sharp$:
\[
B^\sharp:=\{(B_{i,j})_{1\leq i<j\leq n+1}|
B_{i,j}\in \mathbb{Z}\}(=\mathbb{Z}^{\frac{1}{2}n(n+1)}),
\]
Now, we define the action of $\eit^\beta$ $(\beta\geq0)$
on $B^\sharp$ by
\begin{equation}
 \eit^\beta((B_{k,j})_{1\leq k<j\leq n+1})
=((B_{k,j}+\beta_{k,j})_{1\leq k<j\leq n+1}),
\qq\qq
 \beta_{k,j}
:=\begin{cases}
\beta^{(i)}_k&\text{if }j=i,\\
-\beta^{(i)}_k&\text{if }j=i+1,\\
0&\text{otherwise}
\end{cases}
\label{BBB}
\end{equation}
Here note that in the definition of $B^\sharp$,
$B_{i,i}$'s do not appear since the formula $\beta^{(i)}_k$
does not depend on $B_{i,i}$'s as mentioned in the 
remark of the last subsection.

The explicit action of the Kashiwara operator 
$\til e_i$ (resp. $\fit$) 
on $B^\sharp$ is given by (\ref{BBB})
taking $\beta=1$
(resp. $\beta=-1$).
Indeed, the crystal structure of $B^\sharp$ 
is described as follows:

For $v=(B_{i,j})\in B^\sharp$ set 
\begin{eqnarray}
&& b^{(i)}_k(v):=\sum_{1\leq l\leq k}B_{l,i+1}
-\sum_{1\leq l\leq k}B_{l,i},\nn \\
&& \begin{cases}
\displaystyle
\vep_i(v):=\mathop{\rm max}_{1\leq k\leq i}\{ b^{(i)}_k(v)\},\\
\displaystyle
wt(v):=-\sum_{i=1}^{n}(\sum_{1\leq k\leq i
\atop{i+1\leq j\leq n+1}}B_{k,j})\al_i,
\\
\vp_i(v):=\lan h_i,wt(v)\ran+\vep_i(v),
\end{cases}\label{data}
\end{eqnarray}
\begin{eqnarray*}
&&m_i=m_i(v):={\rm min}\{k|1\leq k\leq i, \,\, b^{(i)}_k(v)
=\vep_i(v)\}\\
&&M_i=M_i(v):={\rm max}\{k|1\leq k\leq i, \,\, b^{(i)}_k(v)
=\vep_i(v)\}.
\end{eqnarray*}
The actions of $\til e_i$ and $\til f_i$ on $v=(B_{i,j})$ 
are given by
\begin{eqnarray}
&&\hspace{-50pt}\til f_i:
\begin{cases}
B_{k,j}\longrightarrow B_{k,i}&
\text{if $(k,j)\ne (M_i,i),(M_i,i+1)$}\\
B_{M_i,i}\longrightarrow 
B_{M_i,i}-1&\text{if $(k,j)=(M_i,i)$}\\
B_{M_i,i+1}\rightarrow 
B_{M_i,i+1}+1&\text{if $(k,j)=(M_i,i+1)$}
\end{cases}\label{fityt}\\
&&\hspace{-50pt}\til e_i:
\begin{cases}
B_{k,j}\longrightarrow B_{k,i}
&\text{if $(k,j)\ne (m_i,i),(m_i,i+1)$},\\
B_{m_i,i}\longrightarrow B_{m_i,i}+1&\text{if $(k,j)=(m_i,i)$}\\
B_{m_i,i+1}\rightarrow B_{m_i,i+1}-1&\text{if $(k,j)=(m_i,i+1)$}
\end{cases}\label{eityt}
\end{eqnarray}
\begin{thm4}
By the setting (\ref{data}), (\ref{fityt}) and (\ref{eityt}),
we obtain a free crystal $B^\sharp$.
\end{thm4}
{\sl Proof.}
It suffices to check the axioms (\ref{c1})--(\ref{c5})
in Definition \ref{crystal} and 
the bijectivity of $\eit$ or $\fit$.
Indeed, (\ref{c1})--(\ref{c3}) are trivial from 
(\ref{data}), (\ref{fityt}) and (\ref{eityt}).
The assumption of (\ref{c5}) never occurs. Thus, we may show that 
$\eit\fit={\rm id}=\fit\eit$.
For $v=(B_{i,j})$, set $p:=M_i(v)$, which implies
\[
 b_1^{(i)}(v),\cd, b_{p-1}^{(i)}(v)\leq
 b_p^{(i)}(v)> b_{p+1}^{(i)}(v),\cd, b_i^{(i)}(v).
\]
By the definition of $ b_k^{(i)}$ and the action of $\fit$, 
we have
\[
 b_k^{(i)}(v)=
\begin{cases}
 b_k^{(i)}(v) &1\leq k<p,\\
 b_p^{(i)}(v)+1 &k=p, \\
 b_k^{(i)}(v)+2 &p<k\leq i.
\end{cases}
\]
Thus, we have
\[
 b_1^{(i)}(\fit v),\cd, b_{p-1}^{(i)}(\fit v)<
 b_p^{(i)}(\fit v)\geq b_{p+1}^{(i)}(\fit v),\cd, b_i^{(i)}(\fit v),
\
\]
which means $M_i(v)=p=m_i(\fit(v))$. Similarly, we have 
$m_i(v)=M_i(\eit v)$. It follows from these that 
$\eit\fit={\rm id}_{B^\sharp}=\fit\eit$ 
and then we get (\ref{c4}) and 
the bijectivity of $\eit$ and $\fit$.
\qed

\vskip5mm
\nd
{\sl Remark.} It is unknown whether the crystal graph of $B^\sharp$
is connected or not.
\subsection{Tropicalization of $B^\sharp$}
Let us see that a tropicalization of the crystal $B^\sharp$
is the geometric crystal 
on the unipotent subgroup $U^-\subset SL_{n+1}(\bbC)$
treated in Sect.4.

Applying the following correpondence to (\ref{alik}) and
(\ref{betaik}), and (\ref{gamma-theta}) and (\ref{data})
\begin{eqnarray*}
&&x\cdot y\longleftrightarrow x+y\\
&&x/y\longleftrightarrow x-y\\
&&x+y\longleftrightarrow{\rm max}(x,y)\\
&&i \longleftrightarrow i+1
\end{eqnarray*}
we obtain $\al_{k}^{(i)}\leftrightarrow \beta_k^{(i)}$ and then
\begin{equation}
{\cal UD}_{\hat\theta,T_0}(e_i^c)=\eit^c,\qq
{\cal UD}_{\hat\theta,T_0}(\gamma)=wt,
\label{UD-theta}
\end{equation}
(where $T_0:=(\bbC^\times)^{\frac{n(n+1)}{2}}$),
which implies the following theorem:
\begin{thm4}
We have ${\cal UD}_{\hat\theta,T'}(\chi_{U^-})
=(B^\sharp,wt,\{\eit\}_{i\in I})$, {\it i.e.},
the geometric crystal $\chi_{U^-}$ on 
$U^-\subset SL_{n+1}(\bbC)$ defined in Sect 4 is a tropicalization of 
the crystal $B^\sharp=(B^\sharp,wt,\{\eit\}_{i\in I})$.
\end{thm4}

\end{document}